\documentstyle[11pt, sec, psfig, diagram]{article}

\newtheorem{theorem}{Theorem}[section]
\newtheorem{lemma}[theorem]{Lemma}
\newtheorem{proposition}[theorem]{Proposition}
\newtheorem{corollary}[theorem]{Corollary}
\newtheorem{remark}[theorem]{Remark}
\newtheorem{definition}[theorem]{Definition}
\newtheorem{ex}[theorem]{Example}

\newfont{\Bbb}{msbm10 scaled \magstep1}

\newcommand\uc{\underline{\hbox{\Bbb C}}}

\newcommand\bC{\hbox{\Bbb C}}

\newcommand{\bL}{{\hbox{\Bbb L}}}
\newcommand{\bQ}{{\hbox{\Bbb Q}}}
\newcommand\bR{\hbox{\Bbb R}}
\newcommand\bS{\hbox{\Bbb S}}

\newcommand\bZ{\hbox{\Bbb Z}}

\newcommand\bRP{\hbox{\Bbb RP}}

\newfont{\es}{eusm10 scaled \magstep1}
\newfont{\ses}{eufm8 scaled \magstep1}
\newfont{\gt}{eufb10 scaled \magstep1}
\newfont{\sg}{eufb8 scaled \magstep1}
\newfont{\goth}{eufb10 scaled \magstep2}

\newcommand{\im}{\hbox{\gt Im}}

\newcommand{\gA}{\hbox{\gt A}}
\newcommand{\gB}{\hbox{\gt B}}

\newcommand{\gog}{{\hbox{\gt g}}}

\newcommand{\dir}{\hbox{\es D}}
\newcommand{\en}{\hbox{\es E}}
\newcommand{\uen}{\underline{\hbox{\es E}}}
\newcommand{\hes}{\hbox{\es H}}
\newcommand{\hhes}{\tilde{\hbox{\es H}}}

\newcommand{\Lie}{\hbox{\es L}}
\newcommand{\X}{\hbox{\es X}}

\newcommand{\can}{\hbox{\es K}}

\def\ra{\rightarrow}

\def\be{\begin{equation}}
\def\ee{\end{equation}}

\def\lan{\langle}
\def\ran{\rangle}
\newcommand{\na}{{\nabla}}

\newcommand{\hg}{\hat{g}}
\newcommand{\bc}{{\bf c}}
\newcommand{\hbc}{\hat{\bf c}}
\newcommand{\hsi}{{\hat{\sigma}}}
\newcommand{\hpsi}{\hat{\psi}}
\newcommand{\hA}{{\hat{A}}}

\newcommand{\hconf}{\widehat{\hbox{\gt C}}}
\newcommand{\hgauge}{\widehat{\hbox{\gt G}}}
\newcommand{\conf}{\hbox{\gt C}}
\newcommand{\gauge}{\hbox{\gt G}}
\newcommand{\hmodu}{\widehat{\hbox{\gt M}}}
\newcommand{\modu}{\hbox{\gt M}}

\newcommand{\co}{{\sf C}}

\newcommand{\hco}{\hat{\sf C}}

\newcommand{\hbS}{\hat{\hbox{\Bbb S}}}

\newcommand{\ii}{{\bf i}}

\newcommand{\mod}{ {{\rm mod}\,}}

\newcommand{\dps}{\dot{\psi}}

\newcommand{\dta}{\dot{q}}

\newcommand{\da}{\dot{a}}
\newcommand{\ida}{{\bf i}{\dot{a}}}

\newcommand{\si}{\sigma}

\newcommand{\ve}{{\varepsilon}}
\newcommand{\vfi}{{\varphi}}

\def\dt{\frac{d}{dt}}

\newcommand{\pat}{{\partial_\tau}}

\def\nah{\hat{\nabla}}

\begin{document}

\title{Finite Energy Seiberg-Witten Moduli Spaces on 4-Manifolds Bounding Seifert Fibrations}

\author{Liviu I.Nicolaescu\thanks{{\bf Current address}: Dept.of Math., McMaster University, Hamilton, Ontario,  L8S 4K1, Canada; nicolaes@icarus.math.mcmaster.ca}}

\date{Version 4: November 22, 1997}

\maketitle

\begin{abstract} We compute virtual dimensions of finite energy Seiberg-Witten moduli spaces on $4$-manifolds bounding Seifert fibrations. The key moment is the determination of certain  eta invariants.  As an application of these computations we  compute  the Froyshov invariants of many Brieskorn homology spheres (e.g. $\Sigma(2,3,6k\pm 1)$, $\Sigma(2,4k+1,4k+3)$, $\Sigma(3,3k+1,3k+2)$). In turn, these lead to  interesting topological results. For example, we prove that any   negative definite $4$-manifold bounding a Brieskorn homology sphere $\Sigma(2,3, 6k+1)$  must have diagonalizable intersection form.
\end{abstract}

\addcontentsline{toc}{section}{Introduction}

\begin{center}
{\bf Introduction}
\end{center}

\bigskip

 The goal of the  present paper  is to compute the virtual dimensions of finite energy Seiberg-Witten moduli spaces for $4$-manifolds bounding unions   of Seifert  fibrations.  For cylinders over Seifert manifolds  these  moduli spaces    describe tunnelings between the critical sets of the $3$-dimensional Seiberg-Witten energy functional. These virtual dimensions where computed in [MOY]   by identifying the space of tunnelings with an algebraic-geometric moduli spaces and then using a Riemann-Roch counting argument (see also [KMO] for a spectral flow approach).

We will follow a strategy  similar to the one we used in [N2]  where we dealt with the slightly simpler problem when the $4$-manifold bounds {\em smooth} $S^1$-bundles over  Riemann surfaces.  To compute the virtual dimension we will use the local description in [MMR] of such moduli spaces.   The  computation splits into two main steps.  Firstly,  one needs to understand the  critical sets of the 3-dimensional  Seiberg-Witten functional on the boundary. This was accomplished in [MOY] and [N1] where the critical set was explicitly described  and it was shown  that  (modulo  some very explicit exceptions) they are Bott nondegenerate.     With these informations in hand one can now proceed to compute these dimensions and the results of [MMR]  reduce  this problem to the computation of an Atiyah-Patodi-Singer (A-P-S) index for an operator which near the boundary has the form
\[
{\cal O}= \frac{\partial}{\partial t}- {\hhes}.
\]
Above,  $t$ denotes the outgoing normal coordinate near the boundary. The operator ${\hhes}$ further decomposes as
\[
 {\hhes}={\hhes}_0 +{\cal P}
\]
 where  ${\hhes}_0$ is a direct sum  
\[
{\hhes}_0=({\rm Odd \; signature\; operator}) \oplus ({\rm a \; certain \; Dirac \; operator})
\]
 and ${\cal P}$ is an explicit  zeroth order perturbation of ${\hhes}_0$.  An excision argument reduces the computation of the  A-P-S index  of ${\cal O}$ to   the computation of the A-P-S index of  
\[
{\cal O}_0= \frac{\partial}{\partial t}- {\hhes}_0
\]
and the computation of the spectral flow of the affine path ${\hhes}_0+tP$, $0\leq t \leq 1$. 

The perturbation analysis    employed in [N2] for the computation of $SF({\hhes}_0+t{\cal P})$ extends verbatim to the more general case of Seifert fibration.    The eta invariant of ${\hhes}_0$  is a sum   of the eta invariant of the odd signature operator and the eta  invariant of a  certain Dirac operators.     The eta invariant of the signature operator on Seifert fibrations is essentially a topological invariant  and  was computed in [O].  It can be expressed in terms of the classical Dedekind sums (see [HZ])  determined by the singular fibers  and some topological quantities.   The story is  quite  different for Dirac operators   since  in general  their eta invariants are extremely sensitive to the background geometry.    In  [N2] we discovered that in the   case of smooth circle bundles  the  eta invariants of the Dirac operators relevant in Seiberg-Witten theory  are   very ``blind''.  They detect very little of the background geometry and, moreover, they are insensitive to various couplings.

We could not use the  adiabatic approach in [N2]  since we are not aware of an analogue of Bismut-Cheeger-Dai result for fibrations collapsing onto orbifolds. Instead, we follow a  hands-on approach  outlined in Appendix C of [N2] which relies  on certain algebraic identities satisfied by  the Dirac operators   arising in   Seiberg-Witten theory.  We were even  able to compute  the entire eta functions $\eta(s)$ of these operators.   The singular fibers introduce corrections  to $\eta(0)$ in the form of  some Dedekind-Rademacher sums. (We refer to [Ra] for a presentation of these  objects.)  Surprisingly, the presence of singular fibers restores some of the geometric sensitivity of  the eta invariants in all but one instance namely when the operators are coupled with flat connections.

When specialized to tunnelings, our virtual dimension  formula looks  quite different  from the formula      obtained in [MOY] which involves  some Hirzebr\"{u}ch-Jung continued fractions.    We   did not directly   proved the equality of these two descriptions (although one can conceivably use Rademacher's reciprocity law in [Ra] to achieve this) but numerical  experimentations \footnote[1]{The {\em MAPLE} program used for numerical  experimentations can be downloaded from my homepage at  McMaster University}
 show perfect agreement.  

The eta invariant computations in this paper can be used in conjunction with the adiabatic analysis in [N1]  to provide upper estimates for the Froyshov invariants  $Y_{a,b,c}$ (introduced  in [Fr])   of the Brieskorn homology spheres $\Sigma(a,b,c)$.  In many cases  these estimates are optimal and have interesting topological consequences. for example we prove the following result  (Corollary \ref{cor: liviu}, \ref{cor: liviu1}).
\medskip

\noindent {\em (i) $Y_{2,3,6k+1}=Y_{2,4k+1,4k+3}=Y_{3,3s+1,3s+2}=0$, $Y_{2,3,6k-1}=8$. 

\noindent (ii) $X$ is a   $4$-manifold  with negative definite intersection form  bounding a Brieskorn sphere $\Sigma(2,3, 6k+1)$ then  the intersection form of $X$ must be diagonalizable. (The same result is true for  $\Sigma(2,4k+1,4k+3)$ and $\Sigma(3,3s+1,3s+2)$ and it can be proved directly using  Donaldson's first theorem and the fact proved in [CH] that these Brieskorn spheres bound contractible manifolds.)

\noindent (iii) If $X$ is a $4$-manifold  with negative  definite intersection form $q$ which bounds a $\Sigma(2,3,6k-1)$  and $q$ splits as $q_1\oplus q_2$ with $q_2\neq 0$ even, then $q_2=-E_8$ and $q_1$ is diagonalizable.}

\medskip

The equality $Y_{2,3,5}=8$ is proven in [Fr]  relying on some very special geometric features of the Poincar\'{e} homology sphere $\Sigma(2,3,5)$.

The above result is not entirely  obvious  even when $X$ is the canonical plumbing associated  to $\Sigma(2,3,6k\pm1)$. The  plumbing diagrams are described in Figure \ref{fig: seif1}, \ref{fig: seif2} in \S 3.2.    We discovered   in conversations with Ian Hambleton a simple  combinatorial argument   establishing   these special  facts directly (Remark \ref{rem: funny}).

The paper  is composed of three main sections. Subsection $\S 1.1$ lists the basic topological and geometric facts about Seifert manifolds    used  in the rest of the paper.     In $\S 1.2$ we introduce  the  ``adiabatic'' Dirac operators which as  shown in [MOY] and [N1] play an important role in the $3$-dimensional theory. Most of this section  is occupied with the computation of their eta invariants. We also explain  how to obtain the eta invariants of the traditional Dirac operators (Remark \ref{rem: lcdir}).

The first half of   Section 2  is a brief  survey of the  main results in [MMR] as they apply to our  case.   The virtual dimension formula (\ref{eq: vd5}) is then proved in $\S 2.3$. 

Then third section is devoted to applications. We begin with a brief subsection in which we apply our formul{\ae} to tunnelings.  In the second half, we specialize to Brieskorn homology spheres $\Sigma(a,b,c)$.  Using the adiabatic analysis of  [N1] and the  explicit  knowledge of the eta invariants we produce  ({\em often optimal !}) upper bounds  for the  Froyshov invariants.   These are then applied to obtain  informations about  the intersection forms of $4$-manifolds bounding  such homology spheres. We conclude with a  speculative  subsection   where we formulate several conjectures suggested by numerical  experiments.

\medskip

\noindent{\bf Acknowledgments}   I am indebted to Selman Akbulut and Ron Fintushel for their interest and for generously sharing their expertise with me. Also, I want to thank Nikolai Saveliev for a very helpful e-mail exchange.

\tableofcontents

\section{Differential geometric preliminaries}

\subsection{$V$-bundles and Seifert fibrations}
In this paper, as in [FS] or [MOY], we will think of Seifert  manifolds as sphere bundles  determined by  complex line $V$-bundles  over   Riemann $V$-surfaces.

 The building blocks of $V$-manifolds (or orbifolds) are  quotients of smooth manifolds   a finite group action.  A $V$-manifold is obtained by gluing together such building blocks. More precisely,  according to [Sa] a $V$-manifold   of dimension $n$ is a collection  
\[
M=(\, |M|, (U_i)_{i\in I}, (\tilde{U}_i)_{i\in I}, (G_i)_{i\in I}, (\phi_i)_{i\in I}, (\phi_{ji})\, )
\]
with the following properties.

\noindent $\bullet$ $|M|$ is a Hausdorff space.

\noindent $\bullet$  $(U_i)_{i\in I}$ is an open cover of $|M|$   such that  for any $x\in U_{ij}:=U_i\cap U_j$ there exists $k\in I$  so that $x\in U_k \subset U_{ij}$.

\noindent $\bullet$  $G_i$ denotes  a finite group  acting effectively  on   a connected open set $\tilde{U}_i\subset {\bR}^n$. We assume the fixed point set of $G_i$ has  dimension $\leq n-2$.

\noindent $\bullet$  $\phi_i$ is a $G_i$-invariant map $\tilde{U}_i\ra U_i$ such that the induced map $\tilde{U}_i/G_i\ra U_i$ is a homeomorphisms. 

\noindent $\bullet$  For any pair $(i,j)$ such that $U_i \subset U_j$ there exists  a diffeomorphism $\phi_{ji}$ of $\tilde{U}_i$ onto an open subset of $\tilde{U}_j$ such that the diagram below is commutative
\[
\begin{diagram}
\node{\tilde{U}_i}\arrow{s,l}{\phi_i}\arrow{e,t}{\phi_{ji}}
\node{\tilde{U}_j}\arrow{s,r}{\phi_j} \\
\node{U_i}\arrow{e,t}{\imath}
\node{U_j}
\end{diagram}
\]

\begin{remark}{\rm  The effectiveness of the actions of the groups $G_i$  and the assumption on the fixed point set are usually not included in the definition of  an orbifold.  As shown in Lemmata 1 and 2  of [Sa] these imply that if $U_i\subset U_j \subset U_k$ then  there exists $\gamma_{ki}\in G_k$ such that $\gamma_{ki}\phi_{ki}=\phi_{kj}\phi_{ji}$.}
\label{rem: satake}
\end{remark}

If $N=(|N|, (V_\alpha)_{\alpha\in A}, (\tilde{V}_\alpha), (H_\alpha), (\psi_\alpha), (\psi_{\beta \alpha})\,)$  is another $V$-manifold  then a smooth $V$-map from $M$ to $N$ is a collection $ (\vfi, (h_i))$ with the following properties.

\noindent $\bullet$ $\vfi$ is  a map $I\ra A$ such that  $V_{\vfi(i)} \subset V_{\vfi(j)}$ if and only if  $U_i \subset U_j$.

\noindent $\bullet$  For each $i\in I$ $h_i$ is a smooth map $\tilde{U}_i \ra \tilde{V}_{\vfi(i)}$ such that  if $U_i\subset U_j$ the diagram bellow is commutative  ($\alpha=\vfi(i)$, $\beta=\vfi(j)$.
\[
\begin{diagram}
\node{\tilde{U}_i}\arrow{s,l}{\phi_{ji}}\arrow{e,t}{h_i}
\node{\tilde{V}_{\alpha}}\arrow{s,r}{\psi_{\beta\alpha}} \\
\node{\tilde{U}_j}\arrow{e,t}{h_j}\node{\tilde{V}_\beta}
\end{diagram}
\]
One can check that the collection $(h_i)$ does indeed induce a continuous map $h:|M|\ra |N|$.

 Let $G$ denote a Lie group acting  on a smooth manifold $F$.  A   $V$-bundle   with standard fiber   $F$ and structure group $G$)  consists of two $V$-manifolds 
\[
M =(\, |M|, (U_i)_{i\in I}, (\tilde{U}_i)_{i\in I}, (G_i)_{i\in I}, (\Phi_i)_{i\in I}, (\Phi_{ji})\, )
\]
\[
B=(\, |B|, (V_i)_{i\in I}, (\tilde{V}_i)_{i\in I}, (G_i)_{i\in I}, (\phi_i)_{i\in I}, (\phi_{ij})\, )
\]
and a smooth $V$-map  $\pi=({\bf 1}_I, \pi_i):M \ra B$ with the following properties.

\noindent $\bullet$ For every $i\in I$ the group $G_i$ acts on $F$ as well, its action commutes with the action of $G$  and there exists a $G_i$-equivariant homeomorphism 
\[
\Psi_i : \tilde{U}_i \ra \tilde{V}_i\times F .
\]
$\bullet$ If  $V_i \subset V_j$ there exists a smooth map $g_{ji}: \tilde{V_i} \ra G$  such that 
\[
\tilde{V}_i\times F \ni (v, f)\mapsto\Psi_j\circ \Psi_{ji}\circ \Psi_i^{-1}(v, f)= (v, g_{ji}(v)f)\in \tilde{V}_j\times F.
\]
Moreover, if $V_i\subset V_j\subset V_k$ then there exists $\gamma_{ki}\in G_k$  (as in  Remark \ref{rem: satake}) such that $g_{ki}(\gamma_{ki}\phi_{ki}(p))=g_{kj}(\phi_{ji}(p))g_{ji}(p)$.

When $F$ is a vector space and the actions of $G$ and $G_i$ are all linear  the bundle is called a vector $V$-bundle. All the functorial concepts  related to usual vector bundles (direct sum, tensor product etc.) have a $V$-counterpart. One can also speak of $V$-sections, $V$-metrics and $V$-connections in such vector $V$-bundles. These definitions can be safely left to the reader.  All the above concepts have  holomorphic counterparts defined in an obvious way.   

Most of the standard operations of calculus can be carried out in the context of orbifolds as well.   For example, the Chern-Weil theory has an immediate orbifold extension. We only want to mention here the concept of integration  which as in the smooth case requires an orientability condition.    An orbifold 
\[
M=(\, |M|, (U_i)_{i\in I}, (\tilde{U}_i)_{i\in I}, (G_i)_{i\in I}, (\phi_i)_{i\in I}, (\phi_{ij})\, )
\]
is said to be oriented if  the  sets  $\tilde{U}_i$ are oriented, $G_i$ acts by orientation preserving maps and the gluing maps $\phi_{ji}$ are also orientation preserving.  An $n$-form is a collection $\omega$ of $G_i$-invariant $n$-forms $\tilde{\omega_i}$  such that
\[
  \phi^*_{ji}\tilde{\omega}_j=\tilde{\omega}_i,\;\;\forall U_i\subset U_j.
\]
 Using the Remark \ref{rem: satake}  we deduce
\[
\phi^*_{ki}\tilde{\omega}_k=\phi^*_{ji}\phi^*_{kj}\tilde{\omega}_k
\]
so the above definitions are  non-contradictory.    Now choose a partition of unity $(\alpha_i)$ subordinated to the cover $(U_i)$ and define
\[
\int_M \omega =\sum_{i\in I} \frac{1}{\# G_i}\int_{\tilde{U}_i} \phi_i^*(\alpha_i) \tilde{\omega}_i.
\]
As in the smooth case one can verify the definition is independent of the various choices.

We will be interested  only in compact, oriented 2-dimensional orbifolds without boundary and line $V$-bundles over them.  

If $\Sigma$ is such an orbifold then the isotropy groups $G_i$ can only be cyclic, $G_i \cong {\bZ}_{\alpha_i}$, $\alpha_i >2$. These have only isolated fixed points, $x_1, \ldots, x_n$.  The underlying topological space is a smooth Riemann surface $|\Sigma|$.  Thus we can represent  $\Sigma$ as a collection  $(|\Sigma|, \vec{\alpha},\vec{x})$ where $\vec{\alpha}=(\alpha_1, \ldots, \alpha_n)$ are the orders of isotropy and  $\vec{x}=(x_1,\dots, x_n)$ indicates the location of the singular points.  The underlying topological space is completely described by the genus $g$ of $|\Sigma|$ so this orbifold is completely described by the collection $(g,\vec{\alpha}, \vec{x})$.  We will use the notation $\Sigma(g, n;\vec{\alpha},\vec{x})$ but in general we will drop any information which is either irrelevant or clear from the context.

Suppose $L\ra  \Sigma(g; \vec{\alpha},\vec{x})$   is a line $V$-bundle. There is a new piece of information  related to the singularity $\vec{x}$ namely the  representation $\rho_i$ of ${\bZ}_{\alpha_i}$ on the fiber ${\bC}$.  It can only be of the form $\rho_i(\exp(\frac{2\pi {\bf i}}{\alpha_i}))= \exp(\frac{2\pi\beta_i {\bf i}}{\alpha_i})$ where
\be
0 \leq \beta_i < \alpha_i,\;\; i=1,\ldots ,n.
\label{eq: norm}
\ee
We again  collect this information in the vector $\vec{\beta}$.  It is important to emphasize that the isotropies $\vec{\beta}$ are  {\em normalized} by the condition (\ref{eq: norm}).
We can now associate to $L$ its {\em rational degree}  defined as
\[
\deg L =\int_{|\Sigma|} c_1(L) \in {\bQ}.
\]
As in the smooth case   we have
\[
\deg L_1\otimes L_2 =\deg L_1 + \deg L_2
\]
and that is why we will often use the additive notation for the tensor product of line bundles.

Given an orbifold  $\Sigma(\alpha, x)$ with a single singular point $x$  (of isotropy $\alpha$)  we can construct a line $V$-bundle $H_{\alpha, x}$  with  rational degree $1/\alpha$ as follows.   Let $|\Sigma|_x$ denote the   surface obtained by deleting a small open disk $D_x$ around $x$.  The total space of $H_x$   obtained by gluing $|\Sigma|_x \times {\bC}$ to $ \bar{D}_x\times {\bC}$ via the map 
\[
T_\alpha:S^1\times {\bC}\cong \partial |\Sigma|_x \times {\bC} \ra \partial \bar{D}_x\times {\bC}  \cong S^1 \times {\bC}
\]
where
\[
T_\alpha(\exp (\ii \theta), \rho \exp(\ii \vfi) ) = (\, \exp ({\ii}(-\alpha \theta+\vfi)), \rho\exp({\ii}\theta) \,).
\]
Given a line bundle  $L\ra \Sigma(\vec{\alpha},\vec{x})$ with singularity data $\vec{\beta}$ we can form
\[
|L|=L\otimes H_{x_1}^{-\beta_1}\otimes \cdots \otimes H_{x_n}^{-\beta_n}.
\]
One can show (see [FS]) that $|L|$ is a genuine {\em smooth} line bundle over $|\Sigma|$. Its degree is an integer and  we have the following equality
\be
\deg L =\deg |L| +\sum_i \frac{\beta_i}{\alpha_i}.
\label{eq: deg}
\ee
The rational degree and the singularity data completely determine the topological type of a line $V$-bundle. More precisely we have the following result (see [FS] for more details).

\begin{proposition}{\rm  Denote by ${\rm Pic}^t(\Sigma)$ the  space of isomorphisms classes of line $V$-bundles over the orbifold $\Sigma(\vec{\alpha}, \vec{x})$. Define
\[
\tau=\tau_{\vec{\alpha}}: {\rm Pic}^t(\Sigma) \ra {\bQ}\oplus {\bZ}_{\alpha_1}\oplus \cdots \oplus {\bZ}_{\alpha_n}
\]
 by
\[
L(\vec{\beta}) \mapsto (\deg L, \beta_1 (\mod \alpha_1), \ldots, \beta_n (\mod \alpha_n )\, )
\] 
 and $\delta: {\bQ}\oplus {\bZ}_{\alpha_1}\oplus \cdots \oplus {\bZ}_{\alpha_n}$ by 
\[
(c,  \beta_1 (\mod \alpha_1), \ldots, \beta_n (\mod \alpha_n )\, ) \mapsto \left(c-\sum_i\frac{\beta_i}{\alpha_i}\right) (\mod {\bZ}).
\]
Then 
\[
0
\ra {\rm Pic}^t(\Sigma)\stackrel{\tau}{\ra}  {\bQ}\oplus {\bZ}_{\alpha_1}\oplus \cdots \oplus {\bZ}_{\alpha_n} \stackrel{\delta}{\ra} {\bQ}/{\bZ} \ra 0
\]
is a short exact sequence of   abelian groups.}
\label{prop: picard}
\end{proposition}
The line bundle with singularity data $\vec{\beta}$ and rational degree $c$ will be denoted by $L(c, \vec{\beta})$.

A complex $2$-orbifold $\Sigma(g; \vec{\alpha})$ has an associated canonical  $V$-line bundle $K_\Sigma=\Lambda^{1,0}T^*\Sigma$.   Given a holomorphic line $V$-bundle $L=L(\vec{\beta})$ we can form the  Dolbeault  cohomology $H^*(L)$ and  the  Serre duality continues to hold
\[
H^0(L)\cong H^1(K_\Sigma-L),\;\; H^1(L)\cong H^0(K_\Sigma- L).
\]
The Riemann-Roch theorem has an orbifold version due to Kawasaki [Kaw] ($h_0(L):=\dim H^0(L)$)
\be
h_0(L) -h_0(K_\Sigma-L)= 1-g + \deg |L|=1-g+\deg L-\sum_i \frac{\beta_i}{\alpha_i}.
\label{eq: rr}
\ee
In particular if $L=|L|$ is the trivial line bundle  we deduce from the  Riemann Roch theorem that $h_0(K_\Sigma)=g$. Now, using the Riemann-Roch theorem for  $K_\Sigma$ we deduce
\[
\deg |K_\Sigma|=2g-2.
\]
One can verify easily that the singularity data of $K_\Sigma$ are $(\alpha_1-1, \ldots ,\alpha_n-1)$. Hence
\[
\deg K_\Sigma = 2g-2 + \sum_i\left(1-\frac{1}{\alpha_i}\right).
\]
The rational Euler characteristic of $\Sigma$ is defined  as $\chi(\Sigma)=-\deg K_\Sigma$ from which we conclude
\be
\chi(\Sigma)=2-2g -  \sum_i\left(1-\frac{1}{\alpha_i}\right).
\label{eq: euler}
\ee
 A Seifert manifold is by definition the unit sphere $V$-bundle of a complex line $V$-bundle over a 2-orbifold. More precisely, consider the 2-orbifold $\Sigma(g,\vec{\alpha})$ and the line bundle 
\[
L(c,\vec{\beta})\ra \Sigma
\]
such that  if  $\beta_i\neq 0$  then $g.c.d.(\beta_i, \alpha_i)=1$.  The  unit sphere bundle $S(L)$ of   $L$ is a Seifert manifold we will denote by $N (g;b;\vec{\alpha},\vec{\beta})$ where $b =\deg |L|$.  The collection  $(g;b;\vec{\alpha},\vec{\beta})$ is known as the {\em (normalized) Seifert invariant}.  $\Sigma$ is known as the base of  $N$.

The basic topological invariants of a Seifert manifold are known (see [FS]). 

\begin{theorem} {\rm If $N= N(g;b, \vec{\alpha}, \vec{\beta})\,)$ ($m$ singular points) then $\pi_1(N)$ admits the presentation
\[
\lan a_i,b_i,g_i,f,\;i=1,\ldots, g, j=1,\ldots , m\; ; 
\]
\be[a_i,f]=[b_i,f]=[g_j,f]=g_j^{\alpha_j}f^{-\beta_j}=f^{b}\prod [a_i, b_i]\prod g_j=1\ran.
\label{eq: pi1}
\ee
where $f$ denotes the homotopy class of a regular fiber oriented by the action of $S^1$.  The integral cohomology groups are given by
\[
H^1(N, {\bZ})=\left\{
\begin{array}{rlc}
{\bZ}^{2g}&,& \deg {\bL}_0 \neq 0 \\
{\bZ}^{2g+1}&,& \deg {\bL}_0 =0
\end{array}
\right.
\]
and
\[
H^2(N, {\bZ}) \cong \left({\rm Pic}^t(\Sigma)/{\bZ}[L]\right) \oplus {\bZ}^{2g}.
\]
Moreover,  the projection $\pi:N\ra \Sigma$ pulls back line $V$-bundles on $\Sigma$ to genuine {\em smooth} line bundles on $N$ and the subgroup ${\rm Pic}^t/{\bZ}[L]$ of $H^2(N, {\bZ})$ can be identified with the  image of the morphism}
\[
{\rm Pic}^t(\Sigma)\stackrel{\pi^*}{\ra}\{Smooth\;line\; bundles\; on\; N\}\stackrel{c_1}{\ra}H^2(N,{\bZ}).
\]
\label{th: topo}
\end{theorem}

\subsection{Eta invariants of Dirac operators on Seifert manifolds}
In the sequel we fix a $2$-orbifold $\Sigma=\Sigma(g,m;\vec{\alpha},\vec{x})$  and a line $V$-bundle 
\[
{\bL}_0={\bL}_0(\ell, \vec{\beta})\stackrel{\pi}{\ra} \Sigma.
\]
 such that $0<\beta_i <\alpha_i$ and $g.c.d.(\alpha_i,\beta_i)=1$, $\forall i=1,\ldots , m$ and where $\ell$ the rational degree of ${\bL}_0$.  Denote by  $N$ the associated Seifert manifold $N=S({\bL}_0)$. We orient it as the boundary of a complex  manifold  following the convention
\[
 outer \; normal\; \wedge \; {\bf or}\, (boundary) =  {\bf or} \,( manifold).
\]

As explained in [S], the manifold  $N$ admits a locally homogeneous Riemann metric $g_N$.   The natural $S^1$ action preserves  such a metric.  We denote by $\zeta$ the infinitesimal generator of this action.  $\zeta$ is a nowhere vanishing  Killing vector field. This metric induces  $V$-metric $g_\Sigma$ on the base  which we normalize so that ${\rm vol}\,(\Sigma)=\pi$.  $g_\Sigma$ has constant sectional curvature. By eventually rescaling the metric $g_N$ in the $\zeta$ direction we can assume  $|\zeta|_{g_N}\equiv 1$.    Now denote by $\vfi$ the  {\em global angular form}  defined as the $g_N$-dual of $\zeta$. As shown in [MOY] and  [N1] we have a fundamental  identity
\be
d\vfi =-2\ell \ast\vfi.
\label{eq: str1}
\ee
Denote by $Cl(T^*N)$ the bundle of Clifford algebras generated by  $T^*N, g_N)$.  The bundle  $\Lambda^*T^*N$  is naturally a bundle of $Cl(T^*N)$ modules (see [BGV])  and we denote by 
\[
{\bf c}: Cl(T^*N)\ra {\rm End}\,(\Lambda^*T^*N)
\]
the corresponding Clifford multiplication.  The symbol map 
\[
\sigma: Cl(T^*N)\ra  \Lambda T^*N,  u \mapsto {\bf c}(u)\cdot {\bf 1}
\]
is a bundle  isomorphism with inverse known as the {\em quantization map}.   This allows us to define an action of $\Lambda T^*N$ on itself by
\[
\Lambda T^*N \ra {\rm End}\,(\Lambda^*T^*N),\;\; \omega \mapsto {\bf c}(\sigma^{-1}(\omega)).
\]
For simplicity we continue to denote this map with ${\bf c}$. We call the resulting operation the  Clifford multiplication by a form.

 Let $\lan \vfi \ran^\perp$ denote the orthogonal complement of the real line sub-bundle of $T^*N$ spanned by $\vfi$.   As shown in [N1]  the bundle $\lan \vfi \ran^\perp$ is  ${\bf c}(\ast \vfi)$-invariant.   The bundle $\lan\vfi \ran^\perp$ equipped with the (almost) complex structure  $-{\bf c}(\ast \vfi)$    will be called the {\em canonical line bundle} of $N$ and will be denoted by ${\can}$.   We  have the isomorphism
\[
{\can} \cong \pi^* K_\Sigma.
\]
In [N1], $\S 2.2$ we showed  that  ${\can}$  determines a canonical $spin^c$ structure on  $N$ with associated bundle of spinors
\be
{\bS}_{can} \cong  {\can}^{-1} \oplus {\uc}
\label{eq: canonical}
\ee
where  ${\uc}$ denotes generically the trivial complex line bundle. This allows us to identify  the space  ${\rm Spin}^c(N)$ of $spin^c$ structures on $N$ with the topological Picard group, ${\rm Pic}^t(N)$.    The bundle of spinors corresponding to the line bundle $L\ra N$ is
\be
{\bS}_L={\bS}_{can} \otimes L \cong  L\otimes {\can}^{-1} \oplus L.
\label{eq: str2}
\ee
As above, we get a Clifford multiplication map ${\bf c}:\Lambda^*T^*N\ra {\rm End}\,({\bS}_L)$. We  ``orient'' it using the conventions of [BC]
\[
{\bf c}(dv_N)= -{\bf id}.
\]
On  $TN$ there are two natural $g_N$-compatible connections. The first one is  the Levi-Civita connection $\nabla^{LC}$. The other connection $\nabla^\infty$, called the {\em adiabatic connection} in [N1] and [N2] is described  as follows. Using the decomposition 
\[
T^*N = \lan \vfi \ran \oplus \lan \vfi \ran^\perp \cong  \lan \vfi \ran \oplus \pi^*T\Sigma
\]
we define $\nabla^\infty$ as the direct sum of the trivial connection on $\vfi$ and the pullback of the Levi-Civita connection of $T\Sigma$ on $\lan \vfi \ran$.

The line bundle ${\can}^{-1}$ comes with a natural hermitian connection induced by pullback from the Levi-Civita connection on the base.   Thus  a connection on $\det {\bS}_L\cong L^2\oplus {\can}^{-1}$  can be specified by indicating a connection  on $L$.  Fix such a connection $A$.  Using the connection $\nabla^{LC}$ we obtain  a Dirac operator ${\dir}_A$ on ${\bS}_L$ while the adiabatic connection $\nabla^\infty$ induces a different Dirac operator, $D_A$. We will call $D_A$ the  {\em adiabatic Dirac operator}. These two  Dirac operators are related  by the equality
\be
D_A= {\dir}_A + \frac{\ell}{2}.
\label{eq: str3}
\ee
$D_A$ can be  decomposed as
\[
D_A= Z_A+T_A
\]
where,  with respect to the decomposition  (\ref{eq: str2}), the operators $Z$ and $T$ have the matrix descriptions
\be
Z_{A}=\left[
\begin{array}{cc}
{\bf i}\nabla^{A}_{\zeta} & 0 \\
0  & -{\bf i}\nabla^{A}_{\zeta}
\end{array}
\right]
\label{eq: str4}
\ee
and
\be
T_A=\left[
\begin{array}{cc}
0 & {}^\flat\bar{\partial}_A \\
{}^\flat\bar{\partial}_A^* & 0 
\end{array}
\right].
\label{eq: str5}
\ee
We refer to [N1] for the  exact definitions of $Z$ and $T$. It suffices to say that $Z$ uses only derivatives along the fiber direction while $T$  uses only derivatives along horizontal directions. Moreover $Z$ and $T$ interact in  an especially nice manner:
\be
\{Z, T\}:=ZT+TZ=0.
\label{eq: str6}
\ee
The above  equality is responsible (among many other things) for the following nice result.

\begin{proposition}{\rm   Consider a  line $V$-bundle  $L=L(c,\vec{\gamma}) \ra \Sigma=\Sigma(g,m;\vec{\alpha},\vec{x})$  equipped with  a hermitian  connection $A$.   Denote by  $\tilde{L}$ and resp. $\tilde{A}$ the pullbacks of $L$ and $A$ to $S({\bL}_0)$.  

 If $\ell=\deg {\bL}_0 \neq 0$  then the eta function $\eta_{L,A}(s)$ of the Dirac operator $D_{\tilde{A}}$ on ${\bS}_{\tilde{L}}$ is given by
\be
\eta_{L,A}(s)=-2\ell \zeta(s-1)  +\sum_{i=1}^m\frac{1}{\alpha_i^s}\sum_{r=1}^{\alpha_i-1}\left(\{\frac{\gamma_i +r\beta_i}{\alpha_i} \}-\{\frac{\gamma_i-r\beta_i}{\alpha_i}\} \right)\zeta(s, r/\alpha_i).
\label{eq: etapull}
\ee
The various quantities  which appear in the above formula are  defined as follows. $\zeta(s)$ is  Riemann's zeta function,   $\zeta(s,a)$  ($a>0$) denotes the Riemann-Hurwitz function 
\[
\zeta(s,a)=\sum_{n=0}^{\infty}\frac{1}{(n+a)^s}
\]
 and for each real number  $x$ we have denoted  by $\{x\}$  its fractional part defined by $\{x\}\in [0,1)$, $x-\{x\}\in {\bZ}$. In particular,
\[
\eta_{L,A}(0)=  \frac{\ell}{6} -\sum_{i=1}^m (S_i^+-S_i^-)
\]
where $S_i^\pm$ denotes the   sum}
\[
S_i^\pm= \sum_{r=1}^{\alpha_i} \left\{\frac{\gamma_i\pm r\beta_i}{\alpha_i}\right\}\left(\left(\frac{r}{\alpha_i}\right)\right),\;\;\;((x)):=\left\{
\begin{array}{rcl}
\{x\}-\frac{1}{2} & {\rm if} & x\in {\bR}\setminus {\bZ} \\
0 & {\rm if } & x\in {\bZ}
\end{array}
\right.
\]
\label{prop: eta1}
\end{proposition}

\noindent {\bf Proof} \hspace{.3cm}  Define for $i=1, \ldots, m$ the $\alpha_i$-periodic function $G_i :{\bR} \ra {\bR}$
\[
G_i(t)=-\{ \frac{\gamma_i-t\beta_i}{\alpha_i}\}
\]
and set $G(t)=\sum_i G_i(t)$. Note that for any integer $n$  we have 
\[
G(n)=  \deg |L-n{\bL}_0|  -  \deg (L-n{\bL}_0).
\]
Now set $F_i(n)=G_i(n)-G_i(-n)$,  $F(n)= G(n)-G(-n)$.

The connection  $A$ on  $L$ induces a holomorphic structure on $L$ and we denote by ${\bL}$ the holomorphic line  $V$-bundle thus obtained.   The argument in the Appendix  C of [N2] extends to the case of Seifert manifolds (because of Proposition 5.1.3 of [MOY])
.    We deduce (where  $K=K_\Sigma$)
\[
\eta_{L,A}(s)=
\]
\be
\sum^\infty_{n=1}\frac{(\,h_0(K-{\bL}-n{\bL}_0)-h_0({\bL}+n{\bL}_0)\,)+(\,h_0({\bL} -n{\bL}_0)-h_0(K-{\bL}+n{\bL}_0)\,)}{n^s}
\label{eq: trick}
\ee
({\rm Riemann-Roch-Kawasaki})
\[
=\sum^\infty_{n=1} \frac{-( 1-g+ \deg |L+n{\bL}_0| ) +(1-g +\deg |L-n{\bL}_0|)}{n^s}
\]
\[
= \sum^\infty_{n=1} \frac{\deg |L-n{\bL}_0|- \deg |L+n{\bL}_0|}{n^s}
\]
\[
=  \sum^\infty_{n=1} \frac{\deg (L-n{\bL}_0)- \deg (L+n{\bL}_0) + F(n)}{n^s}
\]
\[
= \sum^\infty_{n=1} \frac{-2n\deg {\bL}_0  + F(n)}{n^s}
\]
\[
= -\sum^\infty_{n=1} \frac{2\ell}{n^{s-1}} +\sum_{i=1}^m \sum_{n=1}^\infty \frac{G_i(n)-G_i(-n)}{n^s}
\]
At this point we use the following elementary lemma whose proof can be safely left to the reader.

\begin{lemma}{\rm Suppose $f:{\bR}\ra {\bR}$ is a periodic function of period $p\in {\bZ}_+$. Then}
\[
\sum_{n=1}^\infty \frac{f(n)}{n^s}= \sum_{r=1}^p \frac{f(r)}{p^s}\zeta(s, r/p).
\] 
\end{lemma}
Using the lemma we deduce
\[
\sum_{n=1}^{\infty}\frac{G_i(n)-G_i(-n)}{n^s}=\sum_{r=1}^{\alpha_i}\frac{G_i(r)-G_i(-r)}{\alpha_i^s}\zeta(s, r/\alpha_i).
\]
The $\alpha_i$-th term in the second sum  vanishes since $G_i(\alpha_i)=G_i(-\alpha_i)=\{\gamma_i/\alpha_i\}$. The first part of the proposition is proved.  The second part  follows from the identities $\zeta(0, a)=1/2 -a=-((a))$  $(0 < a < 1)$,   $\zeta(0)=-1/12$ (see for example Chap. XIII of [WW]).  $\Box$

\medskip

\begin{remark}{\rm  The above result  should be compared with the one in the smooth case (no singular fibers).     Suppose we have {\em different} $V$-line bundles $L_1$, $L_2$  on $\Sigma$  equipped with connection $A_1$, $A_2$ which pullback to the same line bundle $\tilde{L}$ on $N$.  Correspondingly we get two connections $\tilde{A}_1$ and $\tilde{A}_2$ and two adiabatic Dirac operators $D_1$, $D_2$. In the smooth case these operators  have the same eta invariant which may not be  the case  when singular fibers are present.    To put it differently, the eta function is sensible  to the coupling connection.  However,  this dependence is mild if $A'_1$ is another connection on $L_1$ then the Dirac operators $D'_1$ associated to $\tilde{A}'_1$ has the same eta function as $D_1$.   The eta function ``has only a vague idea''  which $V$-line bundle on $\Sigma$ was used to construct the Dirac operator. Note that  Serre duality implies immediately that $\eta_{L}(0)=\eta_{K-L}(0)$.  One can verify this directly   using the explicit description of the eta invariant.

(b) The assumption $\ell \neq 0$  in Proposition \ref{prop: eta1} can be dropped  although the above  approach does not cover the case $\ell =0$. In this case $N$ can be represented as a  quotient $\tilde{N}/{\bZ}_k$ where $\tilde{N}\cong S^1\times N$ and ${\bZ}_k$ acts freely and maps fibers to fibers (Thm. 5.4 in [S]).  The identity (\ref{eq: trick}) can be obtained  as in Appendix C of [N2], working equivariantly on $\tilde{N}$. The details can be safely left to the reader. }
\label{rem: duality}
\end{remark}

 The sum  $S_i^+ -S_i^-$ can be expressed in terms of the Dedekind-Rademacher sums introduced in [Ra].   To see this note denote by $q_i$ the inverse of $\beta_i$ modulo $\alpha_i$ i.e. $\beta_iq_i \equiv 1$ (mod $\alpha_i$).   Then
\[
\left\{ \frac{\gamma_i\pm r\beta_i}{\alpha_i} \right\}=
\left\{
\begin{array}{rcl} 
\left(\left( \frac{\gamma_i\pm r\beta_i}{\alpha_i} \right)\right) +\frac{1}{2} & {\rm if} & r\not\equiv \mp q_i\gamma_i\;(\alpha_i) \\
\left(\left( \frac{\gamma_i\pm r\beta_i}{\alpha_i} \right)\right) &{\rm if} & r\equiv \mp q_i\gamma_i\;(\alpha_i)
\end{array}
\right.
\]
We deduce
\[
S_i^\pm = \sum_{r=1}^{\alpha_i}\left(\left( \frac{\gamma_i\pm r\beta_i}{\alpha_i} \right)\right) \left(\left( \frac{r}{\alpha_i} \right)\right) +\frac{1}{2}\sum_{r=1}^{\alpha_i}\left(\left( \frac{r}{\alpha_i}\right)\right) \pm\frac{1}{2}\left(\left(\frac{q_i\gamma_i}{\alpha_i}\right)\right).
\]
The second sum vanishes since the function $x\mapsto ((x))$ is odd. We deduce
\[
S_i^\pm= \sum_{r=1}^{\alpha_i}\left(\left( \frac{\gamma_i\pm r\beta_i}{\alpha_i} \right)\right) \left(\left(\frac{r}{\alpha_i}\right)\right) \pm\frac{1}{2}\left(\left(\frac{q_i\gamma_i}{\alpha_i}\right)\right).
\]
Using the notations in [Ra] we can rewrite
\[
S_i^\pm = s(\pm\beta_i, \alpha_i ;\gamma_i/\alpha_i , 0) \pm\frac{1}{2}\left(\left(\frac{q_i\gamma_i}{\alpha_i}\right)\right).
\]
where
\[
s(\beta,\alpha ; x,y):=\sum_{r=1}^\alpha\left(\left(x +\beta\frac{r+y}{\alpha}\right)\right)\left(\left(\frac{r+y}{\alpha}\right)\right).
\]
The entries $\alpha$, $\beta$ are coprime integers, $\alpha >0$ and $x, y\in {\bR}$.  Note that the above expression depends only on  $x$, $y$ modulo ${\bZ}$.

Thus  the eta invariant of $D_A$ has the form
\be
\eta_A=\frac{\ell}{6} - 2S(\vec{\beta},\vec{\alpha}; \vec{\gamma}) - d(\vec{\beta},\vec{\alpha};\vec{\gamma})
\label{eq: etadede}
\ee
where
 \be
 2S(\vec{\beta},\vec{\alpha},\vec{\gamma}) = \sum_{i=1}^m (\,s(\beta_i,\alpha_i;\gamma_i/\beta_i,0)-s(-\beta_i,\alpha_i;\gamma_i/\alpha_i,0)\,)
\label{eq: S}
\ee
\[
=2\sum_{i=1}^m s(\beta_i, \alpha_i; \gamma_i/\beta_i, 0).
\]
and
\[
d(\vec{\beta},\vec{\alpha};\vec{\gamma})
=\sum_{i=1}^m\left(\left(\frac{q_i\gamma_i}{\alpha_i}\right)\right).
\]
The reason why we prefer  this alternative description of $\eta$ is  due to the Dedekind-Rademacher reciprocity law.  To formulate it we must distinguish two cases.

\noindent  $\bullet$  Both $x$ and $y$ are integers. 
\be
 s(\beta, \alpha; x,y) +s(\alpha, \beta; y,x) =-\frac{1}{4}+\frac{\alpha^2+\beta^2+1}{12\alpha \beta}.
\label{eq: rec1}
\ee
$\bullet$  $x$ and/or $y$ is not an integer. Then
\be
s(\beta,\alpha; x, y)+s(\alpha, \beta; y,x)=((x))\cdot ((y)) + \frac{\beta^2\psi_2(y) +\psi_2(\beta y+\alpha x) +\alpha^2\psi_2(x)}{2\alpha \beta} 
\label{eq: rec2}
\ee
where $\psi_2(x):= B_2(\{x\})$ and $B_2(z)$ is the second Bernoulli polynomial
\[
B_2(z)= z^2-z+\frac{1}{6}.
\]
 Denote by $R(\beta, \alpha; x,y)$ the right hand side  in the above reciprocity identities. Note that $R(\alpha,\beta; y,x)=R(\beta, \alpha; x,y)$.

The reciprocity law coupled with the identities
\[
s((\beta, 1;x,y)=((\beta y +x))\cdot ((y)),\;\;s(\beta, \alpha;x,y)=s(\beta -m \alpha, \alpha ; x +my,y),\;\;\forall m\in {\bZ}
\]
and the Euclid's's algorithm  provides  an efficient  method of computing Dedekind Rademacher sums.   For example
\[
s(4, 7; 2/7,0)=-s(7,4;0,2/7)+ R(4,7;2/7,0)=-s(3,4;2/7,2/7)+R(4,7;2/7,0)
\]
\[
=s(4,3;2/7, 2/7) -R(3,4;2/7,2/7) + R(4,7;2/7,0)
\]
\[
= s(1,3; 4/7,2/7) - R(3,4;2/7,2/7) + R(4,7;2/7,0)
\]
\[
=-s(3,1; 2/7,4/7) + R(1,3;4/7, 2/7) -R(3,4;2/7,2/7) + R(4,7;2/7,0)
\]
\[
-((3\cdot 4/7 + 2/7))((4/7))+  R(1,3;4/7, 2/7) - R(4,3;2/7,2/7) + R(4,7;2/7,0)
\]
\[
=  R(1,3;4/7, 2/7) -R(4,3;2/7,2/7) + R(4,7;2/7,0)= - \frac{3}{28}.
\]

When dealing with reducible solutions of the $3$-dimensional Seiberg-Witten equations   we  encounter Dirac operators  coupled with flat connections. 

We begin by classifying the topological types of such connections. The space of flat connections (modulo gauge transformations) can be identified via the holonomy representation   with ${\rm Hom}\,(\pi_1(N), S^1)\cong {\rm Hom}\, (H_1(N), S^1)$ where for simplicity we  set $H_*(X)$ (resp. $H^*(X))$ denotes the homology (resp. cohomology of $X$) with coefficients in ${\bZ}$.

Given  a representation  $\rho: H_1(N)\ra S^1$ we obtain a line bundle $L_\rho$ equipped with a flat connection.  We get a map
\[
c_1: {\rm Hom}\,(H_1(N), S^1) \ra {\rm Tors}\,(H^2(N)), \;\; \rho \mapsto c_1(L_\rho).
\]
Clearly   if two representations lie in the same component of ${\rm Hom}\,(H_1(N), S^1)$ then they determine isomorphic line bundles.  Using the presentation  (\ref{eq: pi1}) we    see that the component of ${\bf 1}$ in ${\rm Hom}\,(\pi_1(N), S^1)$ is the  group ${\rm Hom}\,(\pi_1(|\Sigma|), S^1)$.  We thus get  a sequence 
\be
1\hookrightarrow {\rm Hom}\, (\pi_1(|\Sigma|), S^1) \stackrel{\pi^*}{\ra} {\rm Hom}\,(\pi_1(N), S^1) \stackrel{c_1}{\ra} {\rm Tors}\, (H^2(N)) \ra 0
\label{eq: flat1}
\ee
in which $\pi^*$  (induced by the projection $\pi:N\ra |\Sigma|$) is an injection,  ${\rm Range}\,(\pi^*) \subset \ker c_1$ and $c_1$ is onto.

\begin{proposition}{\rm   The sequence (\ref{eq: flat1}) is exact.}
\label{prop: flat1}
\end{proposition}

\noindent {\bf Proof} \hspace{.3cm}  Set $\pi_1:=\pi_1(N)$. Since $c_1$ is a surjection  and the group of components of ${\rm Hom}\,(\pi_1, S^1)$ is finite it suffices to show it has the same cardinality as ${\rm Tors}\,(H^2(N))$.

 For any groups $G$, $A$ ($A$-abelian) denote by $H^*(G,A)$ the cohomology of $G$ with coefficients in the trivial $G$-module $A$.   The  short exact sequence 
\[
0 \ra {\bZ} \ra {\bR} \ra S^1 \ra 1
\]
leads to a long exact sequence in cohomology
\[
  H^1(\pi_1, {\bZ})\stackrel{\imath}{\ra} H^1(\pi_1, {\bR}) \ra H^1(\pi_1, S^1)\stackrel{\delta}{\ra}H^2(\pi_1, {\bZ}) \stackrel{\jmath}{\ra} H^2(\pi_1, {\bR}).
\]
Using the  isomorphisms $H^1(\pi_1, A)\cong {\rm Hom}\, (\pi_1/[\pi_1,\pi_1], A)$ for any abelian group $A$  (see for example [Mac])  and the presentation (\ref{eq: pi1})  we deduce
\[
{\rm coker}\imath = {\rm Hom}\, (\pi_1(|\Sigma|), S^1)
\]
and thus we get a short exact sequence
\[
1\hookrightarrow {\rm Hom}\, (\pi_1(|\Sigma|), S^1) \stackrel{\pi^*}{\ra} {\rm Hom}\,(\pi_1(N), S^1) \ra \ker \jmath \ra 0.
\]
Using the universal coefficient theorem  we deduce that $\ker \imath \cong {\rm Tors}\,(H^2(\pi_1))$.   To complete the proof of the proposition it suffices to show that 
\be
{\rm Tors}\,(H^2(\pi_1))\cong {\rm Tors}\,(H^2(N)).
\label{eq: coho}
\ee

Indeed, of the six Seifert geometries  listed in [S], four of them live on contractible spaces.    In these cases $N$ is a $K(\pi_1,1)$ and we have a stronger,   classical isomorphism $H^2(\pi_1)\cong H^2(N)$ (see [Mac]).    The remaining two Seifert geometries are $S^3$ and $S^2\times {\bR}$.  If the universal cover  is  $S^2\times {\bR}$, then the only compact, oriented Seifert manifold covered by $S^2\times {\bR}$ are  $S^2\times S^1$ and ${\bRP}^3\#{\bRP}^3$  (see [S]) with fundamental groups ${\bZ}$ and respectively ${\bZ}_2\ast {\bZ}_2$. Since $H^2({\bZ}) \cong H^2(S^1)\cong 0$  and $H^2({\bZ}_2\ast {\bZ}_2)={\bZ}_2\oplus {\bZ}_2$ ( Thm. VI.14.2 of [HS]) the  isomorphism (\ref{eq: coho}) is trivially satisfied.

Suppose now that $N$ is covered by $S^3$.  Then $\pi_1$ is     finite so  that  $H_1(\pi_1)$ is pure torsion and, using   Cor. IV. 5.5 in [Mac], we get
\[
H^2(\pi_1)\cong {\rm Hom}\,(\pi_1, S^1) \cong {\rm Hom}\,(H_1(\pi_1), S^1)\cong H_1(\pi_1)\cong H_1(N).
\]
The isomorphism  (\ref{eq: coho}) now follows from the Poincar\'{e} duality on $N$. $\Box$

\bigskip

 Let $\tilde{L}\ra N$ denote a line bundle supporting a flat connection $A$.  In particular,  there exists  a line $V$-bundle $L=L(c,\gamma)\ra \Sigma$  which pulls back to $\tilde{L}$.   $A$ defines a holonomy representation 
\[
{\rm hol}_A: \pi_1(N) \ra S^1
\]
and define $\theta=\theta_A\in [0, 1)$ by the equality
\[
\exp (2\pi {\ii} \theta)= {\rm hol}_A(f)
\]
where $f$ denotes the homotopy class of a regular fiber of $N$. Consider the altered connection $B= A + {\ii}\theta \vfi$. It has trivial holonomy along regular fibers.  By proposition   5.1.3 in [MOY] we conclude that there exists a line $V$-bundle $L'=L'(c';\vec{\gamma}')\ra \Sigma)$ and a  $V$-connection $B'$ on $L'$ such that
\[
(\tilde{L}, B) =\pi^*(L', B').
\]
In particular, $ L(c, \vec{\gamma}) \equiv L'(c', \vec{\gamma}')\; ({\rm mod} \, {\bL}_0)$ in ${\rm Pic}^t(\Sigma)$. Thus 
\[
c'-c\in \ell {\bZ}.
\]
 On the other hand $c'=\deg L'$ can be computed  from
\[
c' =\frac{\ii}{2\pi} \int_\Sigma F_{B'}= \ell \theta.
\]
Thus
\be
\ell \theta -c \in {\ell}{\bZ}.
\label{eq: holo}
\ee 
We distinguish two cases.

\medskip

\noindent {\sf A.} $\ell \neq 0$. In this case  the equality (\ref{eq: holo}) coupled with the restriction $0\leq \theta <1$  {\em uniquely} determines $\theta$  by the equality
\be
\theta =\{c/\ell\}.
\label{eq: holo1}
\ee 

\noindent {\sf B.} $\ell =0$. In this case the the only pure torsion line bundle (which can support flats)  are obtained  via pullback from $V$-line bundles of  zero rational degree. Assuming this is the case, the holonomy $\theta$ can have any value  and in fact,  any flat connection on $\tilde{L}$ is homotopic to one with trivial holonomy along fibers.

The  next proposition summarizes the above considerations.

\begin{proposition}{\rm  (a) If $\ell \neq 0$ the quantity $\theta=\theta(\tilde{L}, A)$ is an invariant of the topological type of $\tilde{L}$ which  we will  denote by $\theta(\tilde{L})$. The correspondence 
\[
{\rm Pic}^t(\Sigma)/{\bZ}[{\bL}_0] \ra S^1
\]
given by $L\mapsto \exp (2\pi {\ii} \theta(L))$ is a group morphism. Moreover,  for any flat connection $A$ on $\tilde{L}$  and any $\theta \cong \theta(L)$ $({\rm mod}\,{\bZ})$ the connection $A+{\ii}\theta\vfi$ is the pullback of a  $V$-connection on the line $V$-bundle
\[
L_\theta=L+(\theta-\frac{c}{\ell}){\bL}_0.
\]
(b) If $\ell =0$ there exist flats on $\tilde{L}$ iff $\tilde{L}$ is the pullback of a  zero degree $V$-line bundle on $\Sigma$. In this case for any flat connection $A$ on $\tilde{L}$   the altered connection $A+{\ii}\theta(A, \tilde{L}) \vfi$ is the pullback  of some flat connection on a  line $V$-bundle on the base.}
\label{prop: holo}
\end{proposition}

Set again $\tilde{L}=\pi^*(L(c;\vec{\gamma})\,)$. The spinor bundle ${\bS}_{\tilde{L}}$ has determinant
\[
\det {\bS}_{\tilde{L}}=\tilde{L}^2\otimes{\can}^{-1}.
\]
${\can}^{-1}$ comes equipped with  a hermitian connection $\tilde{B}_{can}$ pulled back from a constant curvature connection $B_{can}$ on $K_\Sigma$ so that $\tilde{B}_{can}$ has trivial holonomy along fibers.    Thus   

Any connection   $\tilde{B}$ on $\tilde{L}$ induces a connection $\tilde{C}=\tilde{B}^{\otimes 2}\otimes \tilde{B}_{can}$ on $\det {\bS}_{\tilde{L}}$.  We are interested in those connections $\tilde{B}$ such that the corresponding   $\tilde{C}$ is flat. This is equivalent to solving (modulo gauge transformations of $\tilde{L}$) the equation  
\be
F_{\tilde{B}}=\frac{1}{2}F_{\tilde{B}_{can}}.
\label{eq: flat2}
\ee
Denote  by ${\cal R}(\tilde{L})$ the  space of gauge equivalence classes of $\tilde{B}'s$ satisfying (\ref{eq: flat2}).  We again consider separately two cases.

\noindent  {\sf A.} $\ell \neq 0$.  It is very easy to construct a particular solution of   (\ref{eq: flat2}). It suffices to pick a constant curvature connection $A_0$ on  $L(c;\vec{\gamma})$. Denote its pullback to $N$ by $\tilde{A}_0$.   Since ${\rm vol}\,(\Sigma)=\pi$ we deduce
\[
F_{\tilde{A}_0}=\frac{2c}{\ii} \pi^*(dv_\Sigma)= \frac{{\ii}c}{\ell} d\vfi.
\]
It is clear now that the connection
\[
\tilde{B_0}=\tilde{A}_0 + {\ii}\frac{\deg K_\Sigma -2c}{2\ell}\vfi
\]
satisfies (\ref{eq: flat2}).    

If $\tilde{B}\in {\cal R}(\tilde{L})$ is another solution then the connection 
\[
\tilde{B}/\tilde{B_0}:=\tilde{B}\otimes\tilde{B}_0^{\otimes -1}
\]
on the trivial bundle is flat. Thus any solution $\tilde{B}$ of (\ref{eq: flat2}) can be represented   as 
\[
\tilde{B}=\tilde{B_0}\otimes ({\rm flat\; connection\;on\;the\;trivial\;bundle})
\]
\be
=\tilde{A_0}\otimes ({\rm flat\; connection\;on\;the\;trivial\;bundle}) + {\ii}\frac{\deg K_\Sigma -2c}{2\ell}\vfi.
\label{eq: flat3}
\ee
 The above  description  has an ambiguity built-in stemming from a choice of a $V$-line bundle $L$ on $\Sigma$ which  pulls back to $\tilde{L}$.  We can  now be more precise about this choice. 

\begin{definition}{\rm   Let $N\stackrel{\pi}{\ra}\Sigma$ denote a Seifert fibration  over a 2-orbifold $\Sigma$ such that  $\ell=\deg (N) \neq 0$. The {\em canonical representative} of a line bundle $\tilde{L}\ra N \in \pi^*(\, {\rm Pic}^t(\Sigma)\,)$ is the $V$-line bundle $L=L(c;\gamma)$ uniquely defined by the requirements  $\pi^*(L)\cong \tilde{L}$ and  
\[
\rho(\tilde{L}):=\frac{\deg K_\Sigma -2c}{2\ell}\in [0,1).
\]  
A {\em canonical connection} on $\tilde{L}$ is the pullback  of a connection $A_0$ on the canonical representative which has constant curvature}
\[
F_{A_0}={\rm const.} dv_\Sigma.
\]
\label{def: can}
\end{definition}
We can rephrase (\ref{eq: flat3}) by saying  that any solution  of (\ref{eq: flat2}) has the form
\be
{\rm canonical\;connection} \otimes {\rm flat\; on\; the \; trivial \; bundle}+{\ii}\rho(\tilde{L})\vfi.
\label{eq: flat4}
\ee
We immediately see the meaning of $\rho(\tilde{L})$: $\exp(2\pi {\ii} \rho)$ is the holonomy along a nonsingular fiber of  an arbitrary connection $\tilde{B}\in {\cal R}$.

\medskip

\noindent {\sf B.} $\ell =0$   Note first that $c_1({\can})$ may not be a torsion class. However there exist  line $V$-bundle $L$ on $\sigma$  $c_1(\pi^*L^2 \otimes {\can}^{-1})$ is torsion. This is equivalent to
\be
2\deg L =\deg K_\Sigma ,\;\;\;L\in {\rm Pic}^t(\Sigma)
\label{eq: cong1}
\ee
If now $\tilde{L}$ is the pullback of a solution $L$ of (\ref{eq: cong1})    we can argue as in the case $\ell \neq 0$ to see that  once we determine a particular solution of (\ref{eq: flat2}) any other can be obtaining by a twist with a flat connection on the trivial line bundle.  The problem is one of existence. Pick any  constant curvature connection $A_0$. Then its pullback to $\tilde{L}$ will satisfy (\ref{eq: flat2}).

The above structural description of the flats  on a Seifert manifold will enable us to compute   the eta functions of   Dirac operators coupled with such connection.   More precisely, we have the following result.

\begin{proposition}{\rm   Consider  an oriented orbifold $\Sigma(g, m;\vec{\alpha})$ with $m$ singular points and let ${\bL}_0$ denote the line $V$-bundle ${\bL}_0=L(\ell;\vec{\beta})$ of degree $\ell \neq 0$  such that $\alpha_i$, $\beta_i$ are coprime $\forall i$.   Denote by $N$ the   associated unit sphere bundle  equipped with the  geometry described  at the beginning of this section. Let $\tilde{L} \in \pi^*{\rm Pic}^t(\Sigma)$ be  a line bundle on $N$ with torsion $c_1$. Denote by $L=L(c,\vec{\gamma})\ra \Sigma$ its canonical representative. denote another $V$-line bundle  over $\Sigma$ and denote by $\tilde{L}$ its pullback to $N$. Equip it with a connection $\tilde{B} \in {\cal R}$  of the form
\[
\tilde{B}= \tilde{B}_0 + {\ii}\rho \vfi
\]
where $\rho = \rho(\tilde{L})=\frac{\deg K-2c}{2\ell}$ and $\tilde{B_0}$ is the pullback of a constant curvature connection $B_0$ on $L$.  Denote by  $\eta_L(s)=\eta_{\tilde{L},\tilde{B}}(s)$ the eta function of the associated  adiabatic Dirac operator $D_{\tilde{B}}$ on ${\bS}_{\tilde{L}}$.

\noindent (a) If $\rho=0$ then
\[
\eta_L(s)=-2\ell \zeta(s-1)  +\sum_{i=1}^m\frac{1}{\alpha_i^s}\sum_{r=1}^{\alpha_i-1}\left(\{\frac{\gamma_i +r\beta_i}{\alpha_i} \}-\{\frac{\gamma_i-r\beta_i}{\alpha_i}\} \right)\zeta(s, r/\alpha_i).
\]

\noindent (b) If $\rho\in (0,1)$ then
\[
\eta_L(s)= \frac{\deg K -\deg|K|}{2}\left( \zeta(s, \rho) -\zeta(s, 1-\rho)\right)
\]
\be
 -\sum_{i=1}^m  \frac{1}{\alpha_i^s}\sum_{k=0}^{\alpha_i-1}\{\frac{\gamma_i-k\beta_i}{\alpha_i}\}\left(\zeta(s,\{\frac{k+\rho}{\alpha_i}\}) -\zeta(s, 1-\{\frac{k+\rho}{\alpha_i}\}) \right)
\label{eq: hreta}
\ee
\[
-\ell\zeta(s-1, \rho) -\ell \zeta(s-1, 1-\rho)
\]
In particular, for $\rho=0$
\[
\eta_L(0)=\frac{\ell}{6}-2S(\vec{\beta},\vec{\alpha}, \vec{\gamma})-d(\vec{\beta},\vec{\alpha},\vec{\gamma})
\]
while if $0<\rho <1$}
\be
\eta_L(0)=  \frac{\deg K -\deg|K|}{2}(1-2\rho)- \sum_{i=1}^m\sum_{k=0}^{\alpha_i-1}\{\frac{\gamma_i-k\beta_i}{\alpha_i}\}(1-2\{\frac{k+\rho}{\alpha_i}\})
\label{eq: ro}
\ee
\[
-\ell \rho(1-\rho) +\frac{\ell}{6}
\]
\label{prop: eta2}
\end{proposition}

\noindent{\bf Proof} \hspace{.3cm} The case $\rho=0$ is contained in Proposition \ref{prop: eta1} so we will consider only the case of fractional holonomy $0<\rho <1$. We will follow the same arguments and notations as in the Appendix C of [N2] and we deduce
\be
\eta(s)=\sum_{\mu>0}\frac{\dim F_\mu -\dim F_{-\mu}}{\mu^s}
\label{eq: etas1}
\ee
where for any $\mu \in {\bR}$ we denoted by $F_\mu$ the  of  pairs
\[
\alpha\oplus \beta \in C^\infty(\tilde{L}) \oplus C^\infty(\tilde{L}\otimes {\can}^{-1})
\]
satisfying
\be
\ii\nabla^{\tilde{B}}_\zeta \alpha =\mu \alpha, \;\; \bar{\partial}_{\tilde{B}}^*\alpha =0.
\label{eq: c6}
\ee
\be
-\ii\nabla_\zeta^{\tilde{B}}\beta=\mu\beta,\;\;\bar{\partial}_{\tilde{B}}\beta=0.
\label{eq: c7}
\ee
Denote by $f_\mu^-$ (resp. by $f_\mu^+$) the dimension of the space of solutions of (\ref{eq: c6})  (resp. (\ref{eq: c7})).  We can rewrite (\ref{eq: etas1}) as
\be
\eta(s)= \sum_{\mu\neq 0}\frac{{\rm sign}\,(\mu) f^+_\mu}{|\mu|^s} + \sum_{\mu\neq 0}\frac{{\rm sign}\,(\mu) f^-_\mu}{|\mu|^s}.
\label{eq: etas2}
\ee
 We will only show how to determine $f_\mu^+$ since the determination of $f_\mu^-$ is entirely similar.

Set for simplicity $\tilde{B}_\mu^\pm= \tilde{B}\mp \ii \mu\vfi$.   Note first that $\bar{\partial}_{\tilde{B}}\beta =\bar{\partial}_{\tilde{B}_\mu^\pm}\beta$  since the transition $\tilde{B}\ra \tilde{B}_\mu^\pm$ does not alter the derivatives along horizontal directions. On the other hand, the equation $-\ii\nabla_\zeta^{\tilde{B}}\beta=\mu \beta$ can be rewritten as
\[
-\ii\nabla_\zeta^{\tilde{B}_\mu^+}\beta=0.
\]
Thus the equations (\ref{eq: c7}) are equivalent to
\be
-\ii\nabla_\zeta^{\tilde{B}_\mu^+}\beta=0, \;\; \bar{\partial}_{\tilde{B}_\mu^+}\beta=0.
\label{eq: c8}
\ee
If  (\ref{eq: c8}) admits a nontrivial solution $\beta$ then $\beta$ must be $\tilde{B}_\mu^+$-covariant constant along the fibers. This implies that the pair $(\pi^*L, \tilde{B}_\mu^+)$ is the pullback of a pair (line bundle $L_\mu$, connection $B_\mu$ on $L_\mu$) on the base $\Sigma$. Using the equality $\tilde{B}^+_\mu=\tilde{B}_0 +(\rho -\mu){\ii}\vfi$ we can determine the curvature of the connection $B_\mu$  from
\[
F_{B_\mu}=F_{B_0}+2\ell(\mu- \rho)\ii dv_\Sigma.
\]
so that 
\[
\deg L_\mu = \deg L -\ell (\mu -\rho)=c-(\mu-\rho).
\]
Thus $\mu -\rho \in {\bZ}$ since $\pi^*L_\mu=\pi^*L$. If we set $n=\mu-\rho$ we deduce $L_\mu=L-n{\bL}_0$.   The connection $B_\mu$ induces  a holomorphic structure on $L_\mu$ and we can now identify $F^+_\mu$ with the space of holomorphic sections of $L_\mu$.   We have thus proved 
\[
f^+_\mu\neq 0 \Rightarrow \mu =n+\rho,\;n\in {\bZ},\;\;f_\mu^+=h_0(L-n{\bL}_0).
\]
Similarly we deduce
\[
f^-_\mu\neq 0 \Rightarrow \mu =n-\rho,\;n\in {\bZ},\;\;f_\mu^-=h_0(K-L-n{\bL}_0).
\]
Using these informations in (\ref{eq: etas2}) we conclude
\[
\eta(s)=\sum_{n\in{\bZ}}\frac{{\rm sign}\,(n+\rho) h_0(L-n{\bL}_0)}{|n+\rho|^s} + \sum_{n\in{\bZ}}\frac{{\rm sign}\,(n-\rho) h_0(K-L-n{\bL}_0)}{|n-\rho|^s}
\]
\[
=\sum_{n\in{\bZ}}\frac{ {\rm sign}\,(n+\rho)\left(\, h_0(L-n{\bL}_0)-h_0(K-L+n{\bL}_0)\, \right)}{|n+\rho|^s}
\]
(use Riemann-Roch-Kawasaki)
\[
=\sum_{n\in {\bZ}}\frac{ {\rm sign}\,(n+\rho)\left( -\frac{1}{2} \deg |K| +\deg|L-n{\bL}_0|\right)} {|n+\rho|^s}
\]
($c=\deg L$, $G(n)$ as in the proof of Proposition \ref{prop: eta1})
\[
=\sum_{n\in {\bZ}}\frac{ {\rm sign}\,(n+\rho)\left( (\deg K -\deg|K|)/2 ) + G(n) -\frac{1}{2}\deg K +c-n\ell\right)}{ |n+\rho|^s}
\]
( $\ell \rho = \frac{1}{2}\deg K -c$)
\[
=\sum_{n\in {\bZ}} \frac{ {\rm sign}\,(n+\rho)\left(\,(\deg(K -\deg|K|)/2 +G(n)\, \right)}{ |n+\rho|^s} - \ell \sum_{n\in {\bZ}} \frac{ {\rm sign}\, (n+\rho) ( n+ \rho)}{ |n+\rho|^s}
\]
($S_\rho:=\rho +{\bZ}$)
\[
=\frac{1}{2}(\deg K-\deg |K|) \sum_{\mu\in S_\rho} \frac{ {\rm sign}\,(\mu)}{|\mu|^s} +\sum_{i=1}^m\sum_{\mu\in S_\rho} \frac{ {\rm sign}\,(\mu) G_i(\mu-\rho)}{|\mu|^s} - \ell \sum_{\mu \in S_\rho} \frac{1}{ |\mu|^{s-1}}
\]
\[
:={\bf S}_1+{\bf S}_2 +{\bf S}_3.
\]
To compute the first two sums we use the following elementary  result.

\begin{lemma}{\rm   For any periodic function $f:{\bR}\ra {\bR}$   of period $p\in {\bZ}_+$  and any $\rho \in (0,1)$ define
\[ 
\eta (s, f,\rho):=\sum_{\mu\in S_\rho}\frac{ {\rm sign}\,(\mu) f(\mu-\rho)}{|\mu|^s}.
\]
Then}
\[
\eta(s,f,\rho)=\sum_{k=0}^{p-1}\frac{f(k)}{p^s}\left(\,\zeta\left(s, \{\frac{k+\rho}{p}\}\right)-\zeta\left(s, 1-\{\frac{k+\rho}{p}\}\right) \,\right).
\]
\label{lemma:  hreta}
\end{lemma}
\noindent {\bf Proof of the lemma} \hspace{.3cm}  For any $\rho \in {\bR}\setminus {\bZ}$ set
\[
S^\pm_\rho = \{ x\in S_\rho\; ;\; \pm x >0\}.
\]
We have
\[
\eta(s,f,\rho) =\sum_{k=0}^{p-1} \frac{f(k)}{p^s}\sum_{n\in {\bZ}}\frac{ {\rm sign}\,(np+k+\rho)}{|n+(k+\rho)/p|^s} =\sum_{k=0}^{p-1} \frac{f(k)}{p^s}\sum_{\mu \in S_{(k+\rho)/p}} \frac{{\rm sign}\,(\mu)}{|\mu|^s}
\]
\[
 =\sum_{k=0}^{p-1} \frac{f(k)}{p^s}\left(\sum_{\mu \in S^+_{(k+\rho)/p} }\frac{1}{|\mu|^s}\right) - \sum_{k=0}^{p-1} \frac{f(k)}{p^s} \left(\sum_{\mu \in S^-_{(k+\rho)/p} }\frac{1}{|\mu|^s}\right)
\]
\[
=\sum_{k=0}^{p-1}\frac{f(k)}{p^s}\left(\zeta(s,\{\frac{k+\rho}{p}\}) -\zeta(s, 1-\{\frac{k+\rho}{p}\}) \right).\;\;\Box
\]

\medskip

Using Lemma \ref{lemma: hreta} we deduce
\[
{\bf S}_1=\frac{\deg K -\deg|K|}{2}\left( \zeta(s, \rho) -\zeta(s, 1-\rho)\right)
\]
and (since $G_i$ is $\alpha_i$-periodic)
\[
{\bf S_2}= \sum_{i=1}^m\sum_{k=0}^{\alpha_i -1}\frac{G_i(k)}{\alpha_i^s}\left(\zeta(s,\{\frac{k+\rho}{\alpha_i}\}) -\zeta(s, 1-\{\frac{k+\rho}{\alpha_i}\}) \right).
\]
Finally, the third sum  can be written as 
\[
{\bf S}_3= -\ell\left(\, \zeta(s-1, \rho) +\zeta(s-1, 1-\rho)\, \right).
\]
To compute $\eta_L(0)$  we only need to use the classical formul{\ae} (see [WW])
\[
\zeta(0,a)=\frac{1}{2}-a,\;\;\zeta(-1,a)=-\frac{1}{12}+\frac{1}{2}a(1-a)\;\;\forall a>0
\]
which imply $\zeta(0, a)-\zeta(0, 1-a)=1-2a$.   The proposition is proved.  $\Box$

\begin{remark}{\rm  (a) The formula we have just proved  illustrates another  surprising feature of the eta invariant.  The choice of $\vec{\gamma}$ in the above proposition is  uniquely determined by $c_1(\tilde{L})$. This shows the eta invariant of an adiabatic Dirac coupled  with a flat connection is a {\em topological} quantity!

(b) The above results extend to the case $\ell =0$.  In this case  $\rho$ should be defined   as half the ``logarithm'' of  the holonomy of the twisting connection along  a regular fiber.}
\end{remark}

 As in  Proposition \ref{prop: eta1} we can express the  eta invariants in (\ref{eq: ro}) in terms of Dedekind Rademacher sums. To do this observe that
\[
\sum_{k=0}^{\alpha_i-1} \{\frac{\gamma_i-k\beta_i}{\alpha_i}\} (1-2\{\frac{k+\rho}{\alpha_i}\})
\]
$(q_i\beta_i\equiv 1$ mod $\alpha_i )$
\[
=\sum_{k=0}^{\alpha_i-1}\left(\; \left(\left(\frac{\gamma_i-k\beta_i}{\alpha_i}\right)\right) +1/2 \;\right) (1-2\{\frac{k+\rho}{\alpha_i}\}) -\frac{1}{2}(1-2\{\frac{q_i\gamma_i+\rho}{\alpha_i}\})
\]
\[
\sum_{k=0}^{\alpha_i-1}\left(\left(\frac{\gamma_i-k\beta_i}{\alpha_i}\right)\right)(1-2\{\frac{k+\rho}{\alpha_i}\}) +\frac{1}{2}\sum_{k=0}^{\alpha_i-1}(1-2\{\frac{k+\rho}{\alpha_i}\}) - \frac{1}{2}(1-2\{\frac{q_i\gamma_i+\rho}{\alpha_i}\})
\]
\[
= \sum_{k=0}^{\alpha_i-1}\left(\left(\frac{\gamma_i-k\beta_i}{\alpha_i}\right)\right)(1-2\{\frac{k+\rho}{\alpha_i}\}) +\frac{\alpha_i-1}{2}+\{ \frac{q_i\gamma_i+\rho}{\alpha_i} \}-\sum_{k=0}^{\alpha_i-1}\{ \frac{k+\rho}{\alpha_i} \}
\]
The last sum is equal to $\rho+\frac{\alpha_i-1}{2}$ and  thus we deduce
\[
\sum_{k=0}^{\alpha_i-1} \{\frac{\gamma_i-k\beta_i}{\alpha_i}\} (1-2\{\frac{k+\rho}{\alpha_i}\})= -\rho  +\{ \frac{q_i\gamma_i+\rho}{\alpha_i} \} +\sum_{k=0}^{\alpha_i-1}\left(\left(\frac{\gamma_i-k\beta_i}{\alpha_i}\right)\right)(1-2\{\frac{k+\rho}{\alpha_i}\}) 
\]
The last sum (denote it temporarily by $S$) can be expressed in terms of  Dedekind-Rademacher sums.   More precisely, since  $0< \rho <1$ is not an integer then
\[
S= 2\sum_{k=0}^{\alpha_i-1}\left(\left( \frac{k\beta_i-\gamma_i}{\alpha_i}\right)\right) \cdot \left(\left(\frac{k+\rho}{\alpha_i}\right)\right)=2s(\beta_i, \alpha_i; \frac{\gamma_i+\beta_i\rho}{\alpha_i}, -\rho).
\]
For $\rho\in (0,1)$ and $\beta q\equiv 1$ mod $ \alpha $ define
\[
F_\rho(\alpha, \beta,\gamma):=\{ \frac{q\gamma+\rho}{\alpha} \}
\]
Then 
\be
\eta(0)= \frac{1}{2}(\deg K -\deg |K|)(1-2\rho)-\ell\rho(1-\rho) +\frac{\ell}{6}+ m\rho
\label{eq: ro1}
\ee
\[
-2 \sum_{i=1}^m s(\beta_i,\alpha_i ;\frac{\gamma_i+\beta_i\rho}{\alpha_i}, -\rho) -\sum_{i=1}^mF_\rho(\alpha_i,\beta_i, \gamma_i) .
\]
Denote for later use
\[
S_\rho(\vec{\beta}, \vec{\alpha},\vec{\gamma})= \sum_{i=1}^m s(\beta_i,\alpha_i ;\frac{\gamma_i+\beta_i\rho}{\alpha_i}, -\rho)
\]
and
\[
F_\rho(\vec{\alpha},\vec{\beta},\vec{\gamma})=\sum_{i=1}^mF_\rho(\alpha_i,\beta_i, \gamma_i).
\]

\begin{remark}{\rm  (a) We want to point out a delicate fact.  Suppose that $N$ is a Seifert ${\bZ}$-homology sphere equipped as usual with a Thurston geometry.  The base of $N$ is an orbi-sphere $S$. Since $H^1(N)\cong H^2(N)\cong 0$ there is a unique  $spin$ and an unique $spin^c$ structure and thus in this case all the spinor bundles ${\bS}_L$ are isomorphic. However there is a plethora of Dirac operators and it is very easy to  confuse them. We propose below an ``accounting'' method.

Denote by ${\bS}_0$ the spinor bundle associated to the unique $spin$ structure. There are two obvious Dirac operators on ${\bS}_0$. One is ${\dir}_0$ obtained  traditionally using the Levi-Civita connection on $TN$ and the other is $D_0$ obtained using the adiabatic connection. They are related by
\[
D_0={\dir}_0+\frac{\ell}{2},\;\;\ell =\deg N.
\]
By using  different connections on $\det {\bS}_0\cong \underline{\bC}$ we can obtain many other.   

Another fundamental pair of examples of Dirac can be obtained in a similar way. If we think of ${\bS}_{can}$ defined in (\ref{eq: canonical})  then 
\[
\det {\bS}_{can}\cong \det{\bS}_0\otimes {\can}^{-1}.
\]
Using the connection $\nabla^{-can}$ on ${\can}^{-1}$ defined as   the  Levi-Civita connection on $K_S^{-1}$  we obtain the operators ${\dir}_{can}$ and $D_{can}$. 

Also, for ``any'' line bundle $L\ra N$ we can regard
\[
\det {\bS}_L\cong \det{\bS}_{can}\otimes L^2\cong {\can}^{-1}\otimes L^2.
\]
Thus any connection $A$ on $L$ induces a connection $A^{\otimes 2}$ on $L^2$ and in this manner we obtained the operators ${\dir}_A$ and $D_A$ we have been studying so far.      More accurately
\[
{\dir}_A={\dir}_0 \otimes(\nabla^{-can}\otimes A^{\otimes 2}
),\;\;D_A=D_0 \otimes(\nabla^{-can}\otimes A^{\otimes 2}).
\]
In particular, we have
\[
D_0=D_{can}\otimes\nabla^{can}.
\]
We have computed  eta invariants for two types of Dirac operators.

\noindent $\bullet$ {\sf Pullback type}
\[
D_{can}\otimes (\pi^*A)^{\otimes2}
\]
where $A$ is a connection on a line $V$-bundle $L\ra S$.

\medskip
\noindent $\bullet$ {\sf Flat type}
\[
D_{can}\otimes A^{\otimes 2}
\]
where $A$ is a connection on the unique line bundle over $N$ such that 
$\nabla^{-can}\otimes A^{\otimes 2}$ is  flat. In this case there is a unique line bundle on $N$ with a unique (up to gauge equivalence) flat connection.  The above operator is none other than ${\dir}_0$!!!

(b) We can say something about the  eta invariant of ${\dir}_0$ as well.  In our case the metric on $N$ is normalized so that the regular fibers have  length $2\pi$. We deduce as in \S 2.3 of [N2]
\[
\eta(D_0)-\eta({\dir}_0) +\frac{\ell}{6}(\ell^2 -\chi) \in 2{\bZ}
\]
where $\chi=-\deg K_S$. If as in [N1] and [N2] we  deform the metric  on $N$ along the fibers so that they have length $2\pi r$ with $0 < r \ll 1$ then we can be more  precise.   Denote  by $D_0(r)$ and ${\dir}_0(r)$ the new Dirac operators  defined in terms of the new metric. Then
\[
\eta(D_0(r))=\eta(D_0)
\]
and the considerations of \S 2.3 in [N2] lead immediately to the equality
\[
\eta(D_0(r))+2h_{1/2}=\eta ({\dir}_0(r)) - \frac{\ell}{6}(\ell^2r^4-\chi r^2)+\left\{\begin{array}{rlc}
4h_{1/2}&, & \ell >0 \\
0       &, & \ell <0
\end{array}
\right. .
\]
Above $h_{1/2}$ denotes the  dimension of the space of holomorphic sections of a possible nonexistent  holomorphic square root of $K_S$. When such a square root does not exist  we set $h_{1/2}=0$.  Note that even   when    such a square root exists  one  verifies immediately that $\deg |K_S^{1/2}| =-1$  so there cannot be any such holomorphic sections. We conclude
\be
\eta ({\dir}_0(r))= \eta(D_0) + \frac{\ell}{6}(\ell^2r^4-\chi r^2).
\label{eq: lcdir}
\ee  
The invariant $\eta(D_0)$  can be computed  using (\ref{eq: ro1})
The  eta invariant of the signature operator on $N$ equipped with such a metric was computed in [O] (watch out for the Seifert invariant conventions there) and it is
\[
\eta^r_{sign}=-\frac{2\ell}{3}(\ell^2r^4-\chi r^2) + \frac{\ell}{3}-{\rm sign}\,(\ell) - 4S(\vec{\beta}, \vec{\alpha})
\]
where 
\[
S(\vec{\beta},\vec{\alpha}):=\sum_{i}s(\beta_i,\alpha_i;0,0)=S(\vec{\beta},\vec{\alpha};\vec{0})
\]
with $S(\vec{\beta},\vec{\alpha};\vec{\gamma})$ defined  in (\ref{eq: S}). We deduce
\be
4\eta({\dir}_0(r))+\eta^r_{sign}= 4\eta(D_0)+ \frac{\ell}{3}-{\rm sign}\,(\ell) -4S(\vec{\beta}, \vec{\alpha}).
\label{eq: fr1}
\ee
This expression is independent $r$!  The reason we considered it is because it appears in the definition of the invariant introduced by Froyshov in [Fr]. The above expression can be alternatively defined as 
\be
{\bf F}(N):=c_1(L)^2-\sigma(V)-8{\rm ind}_{\bf C}\hat{\dir}_{\bf A}
\label{eq: fr2}
\ee
where $V$ is a an oriented simply connected $4$-manifold with $\partial V=N$, $\sigma(V)$ is the signature of $V$,  $L$ is the determinant line bundle associated to a $spin^c$ structure on $V$ and $\hat{\dir}_A$ is a Dirac operator  of this $spin^c$ structure extending the operator ${\dir}_0$. The above index refers to the Atiyah-Patodi-Singer index. Since $N$ is a  ${\bZ}$-homology sphere we deduce from Rohlin's theorem that $c_1(L)^2-\sigma(V)\in 8{\bZ}$. Hence
\be
 4\eta(D_0)+ \frac{\ell}{3}-{\rm sign}\,(\ell) -4S(\vec{\beta}, \vec{\alpha})\in 8{\bZ}.
\label{eq: cong}
\ee
This happens for every Seifert homology sphere $\Sigma(\vec{\alpha})$ (in the notations of [JN]).}
\label{rem: lcdir}
\end{remark}

\begin{ex}{\rm Let us consider in some detail the special case of the Poincar\'{e} homology sphere. It is the Brieskorn sphere $N=\Sigma(2,3,5)$.   The rational degree of $N$ is $\ell=-1/30$. The Seifert invariants  are described in  [JN]. (Note they use different conventions for the Seifert invariants). They are
\[
\vec{\alpha}=(2,3,5), \;\;\vec{\beta}=(1,2,4).
\]
The degree of the canonical class is 
\[
\deg K_S = -2 +1/2+2/3+4/5=-1/30.
\]
The canonical representative of the unique line bundle on $N$ is the trivial line bundle on $S$ which has degree $0$ and singularity data $\vec{\gamma}=(0,0,0)$. The invariant $\rho$ in this case is $1/2$. Using Proposition \ref{prop: eta2}  we deduce
\[
4\eta(D_0)=539/90.
\]
 A simple computation shows that in this case
\[
\frac{\ell}{3}-{\rm sign}\,(\ell) -4S(\vec{\beta}, \vec{\alpha}) = 181/90.
\]
This shows
\[
{\bf F}(\Sigma(2,3,5))=8.
\]
agreeing with (\ref{eq: cong}).   The Froyshov invariant of the Poincar\'{e} homology sphere is also $8$. In \S 3.2  we will describe a general method of producing upper estimates for the Froyshov invariants which we believe  are actually optimal.}
\label{ex: poincare}
\end{ex}

\section{Finite energy Seiberg-Witten monopoles}
\setcounter{equation}{0}

Throughout this section, a hat over an object  will 
signal (unless otherwise indicated) that  it  is a 4-dimensional geometric
object.  

For example, if $N$ is a 3-manifold then on the tube ${\bR}\times N$ there
exist two exterior  derivatives: the 3-dimensional exterior derivative $d$ along the
slices $\{t\}\times N$ and the 4-dimensional exterior derivative $\hat{d}$ so
that $\hat{d}=dt \wedge \partial_t + d$. If $A(t)$ is a family of connection on
some vector bundle $E\ra N$ then we get a bundle $\hat{E}\ra {\bR}\times N$ 
and  we can think of the path $A(t)$ as a  connection $\hA$ on $\hat{E}$.  We
will denote by $F_{A(t)}$ the curvature of $A(t)$ on the slice $\{t\}\times N$
while $\hat{F}_\hA$ will denote the curvature of $\hA$ on the tube.

\subsection{The 4-dimensional Seiberg-Witten equations}
 Let $\hat{N}$ denote an oriented  4-manifold ({\em not necessarily compact}), equipped with a Riemann metric $\hat{g}$. Denote by $\hat{\ast}$ the Hodge star operator
 defined by the metric $\hg$ and the orientation of $\hat{N}$.  {\em Fix} a connection $\nah$ on $T\hat{N}$ compatible with $\hg$. $\nah$ need not be  the Levi-Civita connection.

  Denote by  $Spin_c(\hat{N})$ the  collection of isomorphism classes of $spin^c$ structures on $\hat{N}$.  For each
  $\hsi \in Spin^c(\hat{N})$ we denote by  $\det \hsi$ the associated line bundle and by
 ${\bS}_\hsi ={\hbS}_\hsi^+\oplus {\hbS}_\hsi^-$ the associated bundle of spinors. 
 Note that $\det \hsi \cong \det \hat {\bS}_\hsi^+$.

 Denote by ${\gA}_\hsi$ the space of hermitian connections on ${\bS}_\si$
 compatible  with both the ${\bZ}_2$-grading and the fixed background
 connection $\nah$. More precisely, $A\in {\gA}_\hsi(\hat{N})$ if for any $\alpha \in
 \Omega^1(N)$, any $X\in {\rm Vect}\,(N)$ and any $\hpsi\in C^\infty({\bS}_\si)$ we have
 \[
 \nabla_X^A({\hbc}(\alpha)\hpsi )=\hbc(\nabla_X\alpha)\hpsi
 +\hbc(\alpha)\nabla^A_X\hpsi
 \]
 where 
 \[
 {\hbc}:T^*\hat{N}\ra {\rm Hom}\,({\hbS}_\hsi^+,
 {\bS}_\hsi^-)
 \]
 denotes the Clifford multiplication.  Any connection on $\det \hsi$ determines
 a connection in ${\gA}_\hsi$ and moreover, once we fix a connection $A_0\in
 {\gA}_\hsi(\hat{N})$,  we can identify ${\gA}_\hsi(\hat{N})$ with
 $\ii\Omega^1(\hat{N})$. To any
 connection $\hA\in {\gA}_\hsi(\hat{N})$ we can associate  the Dirac operator
 \[
 \hat{D}_\hA: \Gamma({\hbS}_\hsi^+)\ra \Gamma({\hbS}_\hsi^-)
 \]
 defined as the composition
 \[
 \Gamma({\hbS}_\hsi) \stackrel{\nah^\hA}{\ra} \Gamma(T^*\hat{N}\otimes
 {\hbS}_\hsi^+)\stackrel{\hbc}{\ra} \Gamma({\hbS}_\hsi^-).
 \]

 There is a natural quadratic map 
 \[
 q: \Gamma({\hbS}_\hsi^+)\ra {\rm End}\, ({\hbS}_\hsi^+),\;\;\hpsi \mapsto
 \tau(\hpsi)
 \]
 defined by
 \[
 q(\hpsi)\hat{\phi}=\lan \hat{\phi}, \hpsi\ran -\frac{1}{2} |\hpsi|^2
 \hat{\phi}.
 \]
 In terms of Dirac's bra-ket notation $\tau(\hpsi)$ can be alternatively
 described as 
 \[
 q(\lan\hpsi|)=|\hpsi\ran \lan \hpsi | -\frac{1}{2}\lan \hpsi
 | \hpsi\ran.
 \]
 Note  that for each $\hpsi$ the endomorphism $\tau(\hpsi)$ is symmetric and
 traceless. 
 
 The quantization map from the exterior algebra to the Clifford algebra extends
 the Clifford multiplication to a map 
 \[
 \hbc : \Lambda^* T^*\hat{N} \ra  {\rm End}\, ({\hbS}_\hsi)
 \]
This map has the property that $\hbc(\omega)$ is a traceless, skew-symmetric
endomorphism of ${\hbS}_\hsi^+$ for any  $\hg$-self-dual real valued 2-form
$\omega$.

The Seiberg-Witten equations  (associated to the $spin^c$ structure $\hsi$) are
equations for a pair $(\hpsi, \hA)$ = (spinor in ${\bS}_\hsi^+$, connection in
${\gA}_\hsi(\hat{N})$).  More precisely, they are
\[
(\widehat{SW})\;\;\left\{ 
\begin{array}{rcl}
\hat{D}_\hA \hpsi & = & 0 \\
\frac{1}{2}\hbc(\hat{F}^+_\hA)& = & \tau(\hpsi) 
\end{array}
\right.
\]

In the remaining part of this subsection we will make further  additional
assumptions on the geometry and the topology of $\hat{N}$ and explain how this
affects the Seiberg-Witten equations.

More precisely, assume the manifold $\hat{N}$ can be decomposed as
\[
\hat{N}=\hat{N}_0 \cup [0,\infty) \times N
\]
where $\hat{N}_0$ is a compact oriented 4-manifold with boundary $\partial
\hat{N}_0 =N$.   We will denote  by $t$ the (outgoing) longitudinal coordinate on the cylindrical part of $N$.

Fix a tubular neighborhood  $(-1,0]\times N$ of $N$ in $\hat{N}_0$, a metric
$g$ on $N$ and a connection $\nabla$ compatible with $g$, {\em not necessarily the Levi-Civita connection} of $g$. We assume that along the  infinite cylinder $(-1,\infty)\times N$ the metric $\hg$ is a  product metric $\hg =dt^2 + g$. We fix a   connection $\nah$ compatible with $\hg$ such that along the above cylindrical end  it has the form
 \[
 \nah = \pat \wedge dt + \na .
 \]
 We denoted by $\pat$ the  $\hg$-gradient of $\tau$ where $\tau :\hat{N} \ra [0,\infty)$ is a smooth function
  which coincides with the  canonical projection $[0,\infty)\times N \ra
  [0,\infty)$ on the infinite neck.

  Note that the  $spin^c$ structure $\hsi$ induces  a  $spin^c$
  structure $\si$ on  $N=\partial \hat{N}_0$.  Denote by ${\bS}_\si\ra N$ the
  associated  bundle of spinors and by $\bc :T^*N\ra {\rm End}\, (N)$ the
  corresponding Clifford multiplication. As in the $4$-dimensional case we can
  define ${\gA}_\si(N)$.
  
  Fix a reference connection $\hA_0\in {\gA}_\hsi(\hat{N})$ which  along the neck is
  gauge equivalent to a product connection $dt\otimes \partial_t +A_0$, $A_0\in
  {\gA}_\si(N)$.   Now define the configuration space $\hconf$ as the set of pairs 
  $(\psi,  \hA_0+\ii \hat{a}):=(\hpsi, \hA)$=(spinor, connection) such that    
  $(\hpsi, \ii \hat{a}) \in L^{2,2}_{loc}(\hat{\bS}\oplus{\bf  i}T^*\hat{N})$ and  
  \[
\nah^{\hA}_\pat\hpsi \oplus i_\pat\left(\ii\hat{F}_{\hA}\right) \in
  L^2(\hat{\bS}_\hsi \oplus  \ii\Lambda^1T^*\hat{N}).
  \] 
 We denoted    by $i_\pat$  the contraction by $\pat$.  For brevity,  will denote the elements of $\hconf$ by   the generic symbol $\hco$. 

  \begin{definition}{\rm  (a) A finite energy solution of $(\widehat{SW}_\omega)$
  is a solution $(\hpsi, \hA)$ such that $(\hpsi, \hA-\hA_0)\in \hconf$. 
  
  \noindent (b) A Seiberg-Witten tunneling is a finite energy solution on
  $\hat{N}={\bR}\times N$.}
  \end{definition}
  
  There is an infinite dimensional group $\hgauge$ acting on the configuration space, more precisely 
  \[
  \hgauge =\{ \gamma \in  {\rm Map}\,(\hat{N}, S^1)\; ; \;  \gamma \in
  L^{3,2}_{loc}\}  
  \]
 The group $\hgauge$ acts (on the right) on $\hconf$ and  transforms finite
  energy solutions  to finite energy solutions. Define
  \[
  \hmodu :=\{ (\hpsi, \hA) \; {\rm finite\; energy\; solutions\; of\;}
  \widehat{SW}\}/\hgauge.
  \]
  In this section we want to analyze the the Fredholm properties of the deformation complex naturally associated to 
  $\hmodu$ when $N$ is a Seifert fibration.
  
  We conclude this subsection  with  a simple but crucial observation which
  reveals the dynamical feature of the  Seiberg-Witten equations on cylinders.

   Note that if we set $J=\hbc(d\tau)$
  then  $J$ induces  isomorphisms
  \be
  \hat{\bS}^+_\hsi\!\mid_N\cong \hat{\bS}_\hsi^-\!\mid_N \cong {\bS}_\si
  \label{eq: spin1}
  \ee
  and
  \be
  \bc(\alpha)=J\hbc(\alpha),\;\; \forall\alpha \in \Omega^1(N) \hookrightarrow \Omega^1 ([0,\infty)\times
  N).
  \label{eq: spin2}
  \ee
 A connection $\hA\in {\gA}_\hsi(\hat{N})$ is said to be in a {\em temporal gauge} if
 $i_\pat (\hA- \hA_0)=0$ along  the infinite neck $[0,\infty) \times N$.

 Assume now that $(\hpsi, \hA)$   is a finite energy solution of
 $(\widehat{SW})$ such that $\hA$ is in a temporal gauge.  Along the neck we
 can write
 \[
 \hpsi =\psi(t),\;\; \hA = A_0 + \ii a(t)
 \]
 where $A_0=\hA_0\!\mid_N$, $\psi(t) \in \Gamma ({\bS}_\si)$, $a(t) \in  \Omega^1(N)$, $\forall t\geq 0$.  Then  (along the neck)
 \be
 \hat{F}^+_{\hA} =\frac{1}{2}\{ (F_a +\ast \ii \dot{a}) +dt \wedge (\ii\dot{a} +\ast
 F_a)\}
 \label{eq: sd}
 \ee
 where $A_0 + a(t)$ is the connection on  the line bundle $\det \si$ restricted to  the slice
 $\{t\}\times N$,  $F_a =F_{A_0 +\ii a}$ denotes is curvature and $\ast$ denotes the Hodge star operator on $N$. $A_0+\ii a(t)$  induces a  Dirac operator 
 \[
 {D}_a={D}_{a(t)}: \Gamma({\bS}_\si)\ra \Gamma({\bS}_\si).
 \]
 Using  (\ref{eq: spin1}) and (\ref{eq: spin2}) we deduce that along the neck
 \[
 \hat{D}_\hA= J\left(\partial_t -D_a\right).
 \]
 The equality (\ref{eq: sd})   now implies
 \[
 \hbc(\hat{F}^+_\hA)= \bc (\ast F_a + \ii \dot{a}).
 \]
 Consequently, along the neck, in a temporal gauge,  the Seiberg-Witten equations
 can be rewritten as
 \be
 \left\{\begin{array}{rcl}
 \dot{\psi} &=& {D}_a \psi \\
 \ii\bc (\dot{a})&=&q(\psi)-\frac{1}{2}\bc(\ast F_a) 
 \end{array}
 \right. .
 \label{eq: swtube}
 \ee
 The right-hand-side  of (\ref{eq: swtube})    arises  when one considers the three dimensional counterpart of the Seiberg-Witten equations.

\subsection{The 3-dimensional Seiberg-Witten equations}
To formulate these equations we need to consider a new configuration space. Fix
a connection $A_0\in {\gA}_\si(N)$ and define
\[
\conf =\{(\psi, A)\; ;\; (\psi, (A-A_0) \in L^{1,2}({\bS}_\si \oplus {\bf i}T^*N)\}.
\]
  For brevity, its elements will be denoted by  the symbol $\co$ and we will often write
  $\co=(\psi, a)$ instead of $(\psi, A_0+\ii a)$ whenever no confusion is possible. There is an energy functional  ${\en} :\conf \ra {\bR}$ defined by
\be
{\en}(\psi, A)=\frac{\ii }{4}\int_N a\wedge (F_{A_0}+F_A) +\frac{1}{2}\int_N\lan \psi , D_A\psi\ran dv_g.
\label{eq: energy}
\ee
The gauge group 
\[
\gauge =\{\gamma\in {\rm Map}\,(N, S^1)\; ;\; \gamma \in L^{2,2}\}
\]
acts on $\conf$ and moreover
\[
{\en}(\gamma^{-1}\cdot(\psi,A))-{\en}(\psi, A)=-\int_N\gamma^{-1}d\gamma \wedge F_{A_0} =2\pi{\bf
i}\int_N\gamma^{-1}d\gamma \wedge c_1(A_0)
\]
where  we denoted by $c_1(A_0)$ the  2-from representing the first  Chern class of $\det \si$  associated to $A_0$ via the Chern-Weil construction.  The  $L^2$-gradient  of ${\en}$  is (see [N1] or [MOY])
\[
\nabla {\en}(\psi, A) =\left[
\begin{array}{c}
D_A\psi \\
q(\psi)-\frac{1}{2}\ast F_A
\end{array}
\right]
\]
where  we  tacitly  identified $q(\psi)$ with a purely imaginary  1-form via the Clifford multiplication.   The 3-dimensional Seiberg-Witten equations  can now be  described  as
\[
\nabla{\en} (\co)=0 \Longleftrightarrow \left\{
\begin{array}{rcl}
{\dir}_A\psi & =& 0 \\
{\bc}(\ast F_A) & =& q(\psi)
\end{array}
\right.
\]
We see that  (\ref{eq: swtube}) can be rewritten as a gradient flow equation 
\be
\dot{\co}=\nabla {\en}(\co).
\label{eq: swtube1}
\ee
This last equation suggests that  as $t\ra \infty$  $\co(t)$ converges to a critical point  of ${\en}$.   Assuming the finite energy  condition this  can be proved for arbitrary $N$ using the techniques of [MMR].  However, unlike  the Yang-Mills situation,  the nature of critical points  and  the manner in which they are organized are less transparent in the Seiberg-Witten case.  This is  the reason
 why   we will concentrate on a  special case.

\medskip

{\em In remaining  of the section $N$ will be assumed to be a  Seifert fibration  determined by the $V$-line bundle ${\bL}_0=L(\ell, \vec{\beta})$ of rational  degree $\ell \neq 0$  over the $V$-surface  $\Sigma(g,\vec{\alpha)}$,  $S^1\hookrightarrow N \stackrel{\pi}{\ra} \Sigma$ equipped with  the metric described in $\S 1.2$  As background $g$-compatible connection on $N$ we choose the adiabatic connection $\nabla^\infty$.}

\medskip

 Fix a line bundle $L\mapsto N$  which, as explained in $\S 1.2$, determines a $spin^c$ structure with associated bundle of spinors ${\bS}_L$. Using the decomposition ${\bS}_L\cong \left({\can}^{-1}\otimes L\right)\, \oplus
 \,  L$ we can  represent any section $\psi$ of ${\bS}_L$  as
 $\psi=\psi_- \oplus \psi_+$.  In terms of a local, oriented orthonormal frame $(\zeta,\zeta_1,\zeta_2)$ with dual coframe $(\vfi,\vfi^1,\vfi^2)$  the Seiberg-Witten equations can be  rephrased as  (see [N1] or [MOY])
\be
\left\{
\begin{array}{rrrcl}
{\bf i}{\nabla}^A_\zeta \psi_- & +  \bar{\partial}_A \psi_+ & +\lambda  \psi_- &  = & 0 \\
 & & & & \\
(\bar{\partial}_A)^* \psi_- &  - {\bf i}{\nabla}_\zeta \psi_+  & +\lambda  \psi_+ & = &0 \\
&&&&\\

& &  \frac{\ii}{2}(|\psi_+|^2-|\psi_-|^2) & = &F_A(\zeta_1, \zeta_2) -\frac{1}{2}F_{can}(\zeta_1,\zeta_2)\\
& & & & \\
& & {\bf i}\psi_- \bar{\psi}_+&= &\bar{\ve}\otimes F_A(\zeta_1
+{\bf i} \zeta_2, \zeta)
\end{array}
\right.
\label{eq: sw}
\ee
where ${\ve} =2^{-1/2}(\vfi^1+\ii\vfi^2)$, and $F_{can}$ denotes the pullback of the curvature of $K_\Sigma$ equipped with its Levi-Civita connection.

Set
\[
\conf^*=\{(\psi, A)\in \conf\; ;\; \psi \not\equiv 0\}.
\]
The configurations in  $\conf^*$ are called {\em irreducible}. As in [M]  one  
can show that ${\gB}:=\conf/\gauge$ is a metric space and, moreover,  ${\gB}^* =\conf^*/\gauge$ is  a Banach manifold.  
This is proved  using the  existence of local slices for the $\gauge$-action exactly as in the Yang-Mills case.    For every configuration $\co\in {\conf}$ we will denote by $[\co]$ its image in ${\gB}$.

The solutions of (\ref{eq: sw}) are explicitly described in  [MOY]. 
Here are the relevant facts.

\medskip

\noindent {\bf Fact 1.}  If $c_1(L)$ is not torsion then  (\ref{eq: sw}) has no
solutions.

\medskip

Assume now that $c_1(L) = \kappa \in {\rm Pic}^t(\Sigma)/{\bZ}[{\bL}_0] $  and define
\[
\tilde{R}_\kappa=\{E\in {\rm Pic}^t(\Sigma) \; ; \; 0 < |\deg E -\frac{1}{2} \deg K_\Sigma | \leq \frac{1}{2} \deg K_\Sigma, \;\; E\equiv \kappa \; {\rm mod}\; {\bZ}[{\bL}_0] \}.
\]
For each $E\in \tilde{R}_\kappa$ set $\nu(E)=\deg E -\frac{1}{2}\deg K_\Sigma$. Note that since $\ell \neq 0$ the map $\nu: \tilde{R}_\kappa \ra {\bQ}$ is injective. We will often identify $\tilde{R}_\kappa$ with its image  in ${\bQ}$ via $\nu$.  Now set
\[
R^-_{\kappa}=\{E\in \tilde{R}_\kappa\; ;\; \nu(E) <0,\;\;\deg|E| \geq 0\}
\]
\[
R^+_\kappa=\{E\in \tilde{R}_\kappa\; ;\; \nu(E)>0,\;\;\deg|K_\Sigma -E| \geq 0\}
\]
and
\[
R_\kappa =R^-_\kappa \cup R^+_\kappa.
\]

\medskip

\noindent{\bf Fact 2.}  Any irreducible solution $(\phi, A)$   of (\ref{eq:
sw})  is gauge equivalent  to the pullback of a pair $(\tilde{\phi}, B)$
where $B$ is a connection in a $V$-line bundle $E \ra \Sigma$ in $R_\kappa$ such that $\tilde{\phi}=\tilde{\phi}_-\oplus \tilde{\phi}_-$ is a
section of $C^\infty(K^{-1}\otimes E \, \oplus \, E)$.
The connection $B$ defines holomorphic structures in $K^{-1}\otimes E$ and $E$.
$\tilde{\phi}_-$ is an antiholomorphic section of $K^{-1}\otimes E$  while
$\tilde{\phi}_+$ is a holomorphic section of $E$.    Moreover,  exactly one of
$\tilde{\phi}_-$ or $\tilde{\phi}_+$ is zero according to the identity:
\[
\frac{1}{4\pi}\int_\Sigma(|\tilde{\phi}_-|^2 -|\tilde{\phi}_+|^2) \, dv =\nu(E) \neq 0.
\]
Thus $\tilde{\phi}_+ =0$ iff $\nu(E)>0$ and $\tilde{\phi}_- =0$ iff $\nu(E)<0$. Pairs $(\tilde{\phi}_-\oplus \tilde{\phi}_+, B)$ as above are known as {\em vortex pairs} on $\Sigma$  corresponding to the $V$-line bundle $E$.  If $\nu(E) <0$ we  say we have a {\em holomorphic vortex on} $E$ while if $\nu(E)>0$ we say we have an {\em antiholomorphic vortex} on $E$.

The irreducible  part (mod $\gauge$), denoted by ${\modu}^*$  consists of $\# R_\kappa$ components
\[
{\modu}^*=\bigcup_{E\in R_\kappa} {\modu}_{\kappa, E} =\bigcup_{n\in \nu(R_\kappa)}{\modu}_{\kappa, n}.
\]
  The component ${\modu}_n={\modu}_{\kappa,n}$ corresponding to a choice $\nu(E) =n <0 $ is diffeomorphic to a symmetric product
of $\deg | E|$ copies of $\Sigma$. If $n=\nu(E)>0$  the moduli space is isomorphic to a symmetric product of  $\deg |K -E|$ copies of $\Sigma$. Each component is Bott nondegenerate as a critical set. 

\medskip

\noindent {\bf Fact 3.}  The reducible solutions  consist of pairs  $(0, A)$ where $A$ is a connection on $E$ satisfying (\ref{eq: flat2}).  Modulo $\gauge$ they form a space ${\modu}_\kappa^0$ homeomorphic to the jacobian $J(|\Sigma|)$ which is a $2g$-dimensional torus.   If there exist degenerate  reducible (i.e.  $\ker D_A \neq 0$)    then invariant $\rho$ of  $L\ra N$ must be zero. In this case  the  degenerate solutions can be identified with  the Brill-Noether locus consisting of all  holomorphic  line $V$-bundles of degree $\deg K_\Sigma/2$ which admit nontrivial  holomorphic sections.

\medskip

\noindent{\bf Fact 4.} If  in the definition  of ${\en}$ we fix the reference connection such that ${\en}\equiv 0$   on the reducible part of $\modu$ then the energy along  the component $R_{\kappa, E}$  can be expressed as $const.\nu(E)^2/\ell$, where $const$ is a certain positive universal constant irrelevant  in the sequel.

\bigskip

Associated to each component ${\modu}$ there is a deformation theory
which we now proceed to describe.  We will concentrate only on the irreducible
part $\conf^*$. Since $c_1(L)$ is torsion the energy
functional ${\en}$ is gauge invariant  and  thus it descends to a well defined
functional 
\[
{\uen}: {\gB}^*\ra {\bR}.
\]
The group $\gauge$ is a Hilbert-Lie group and its Lie algebra can be identified
with the space ${\gog}:=L^{2,2}(N, \ii{\bR})$. The exponential map has the form
\[
{\gog}\ni \ii f \mapsto (\exp(\ii f): N \ra S^1).
\]
The tangent space to the orbit ${\cal O}_{\phi, A}$ through  $\co=(\phi, A)$ of the {\em right} action of
$\gauge$  is the range of the infinitesimal action
operator
\[
{\Lie}={\Lie}_{\co}:{\gog}\ra {\X}:=L^{1,2}({\bS}_L)\oplus L^{1,2}(\ii
T^*N),\;\;\ii f\mapsto -\ii f \oplus \ii df.
\]
The tangent space to ${\gB}^*$ at  $[\co]$ can be identified with the orthogonal complement to the tangent space to the orbit ${\cal O}_\co$ and ultimately with the kernel of ${\Lie}^*_\co$, the adjoint of  ${\Lie}_\co$. An integration by parts shows
\[
{\Lie}^*(\dps\oplus \ida)=-\ii d^* a+\ii \im \lan \phi, \dps\ran, \;\;\forall \dps\oplus \ida \in {\X}.
\]
We can use the affine structure of ${\conf}$ to linearize  $\nabla {\en}$ at a given configuration $\co=(\phi, A)$ and we  obtain the {\em unrestricted hessian} at $\co$
\[
{\hes}_\co\left[
\begin{array}{c}
\dps \\
\ida
\end{array}
\right] =\frac{d}{dt}\!\mid_{t=0}\nabla{\en}(\psi +t\dps, A+t\ida)= \left[
\begin{array}{c}
D_A{\dps} +{\bf c}({\ida})\phi  \\
-{\bf i}\ast d{\da} +{\dta}(\phi, {\dps}) 
\end{array}
\right]
\]
The term $ {\dta}(\phi, {\dps}) $   is formally defined by the equality
\[
{\dta}(\phi, {\dps})  :={\dt}\!\mid_{t=0}q(\phi +t{\dps})
\]
where we regard $q$ as a quadratic map $q:{\bS}_L \ra {\bf i}T^*N$.   

The {\em stabilized hessian} of ${\en}$ at $\co=(\phi, A)$ is the unbounded operator on $L^2({\bS}_L\oplus\ii(\Lambda^1\oplus\Lambda^0)T^*N)$ defined by
\[
{\hhes}_{\co}\left[
\begin{array}{c}
{\dps}\oplus{\ida} \\
{\bf i}f 
\end{array}
\right ]:=\left[\begin {array}{cc}
{\hes} & {\Lie} \\
{\Lie}^* & 0
\end{array}
\right]  \left[
\begin{array}{c}
{\dps}\oplus {\ida} \\
{\bf i}f 
\end{array}
\right ] = \left[
\begin{array}{clcl}
D_A\dps  &                               &+ &  {\bc}(\ida)\phi-\ii f\phi   \\
                 &-\ii\ast  d\da +\ii df & +&\dta(\phi, \dps) \\
                  &\ii d^*\da& + & \ii \im \lan \phi, \dps\ran  
\end{array}
\right]
\]
In   [MOY]  it is shown that if $[\co]\in {\modu}_{\kappa, n}$ then  the kernel of the stabilized hessian ${\hhes}_\co$  is naturally isomorphic to the tangent space $T_{[\co]}{\modu}_{\kappa,n}$.  Now define 
\[
{\hhes}_0  \left[
\begin{array}{c}
{\dps}
\\
 {\ida} \\
{\bf i}f 
\end{array}
\right ] =\left[
\begin{array}{crc}
D_A & 0 & 0\\
0 & -\ast d & d \\
0 & d^*  & 0 
\end{array}
\right]  \left[
\begin{array}{c}
{\dps}
\\
 {\ida} \\
{\bf i}f 
\end{array}
\right ] 
\]
and ${\cal P}={\cal P}_\phi$ by
\[
{\cal P}_\phi \left[
\begin{array}{c}
{\dps}
\\
 {\ida} \\
{\bf i}f 
\end{array}
\right ] 
=\left[
\begin{array}{c}
{\bc}(\ida)\phi -\ii f \phi \\
\dta(\phi, \dps) \\
\ii \im \lan \phi, \dps\ran
\end{array}
\right]
\]
Note that  ${\hhes}_0={\hhes}_0(\co)$ is an elliptic selfadjoint operator  for any $\co \in {\conf}$ and ${\hhes}_\co={\hhes}_0+{\cal
P}_\phi$.     For every $\co \in {\conf}$ define $SF(\co)\in {\bZ}$ as the spectral flow of the path ${\hhes}_0(\co) +t{\cal
P}_\phi$, $t\in [0,1]$. 

This spectral flow can be computed exactly as in Sec. $\S 3.3$ of [N2]. More precisely, suppose $\co \in {\modu}_{\kappa, n}$. Thus $\co$ is the pull back of a vortex pair $(\phi_-\oplus \phi_+, B)$   corresponding to a $V$-line bundle $L_\Sigma \ra \Sigma$ in $R_\kappa$ with $\nu(L_\Sigma)=n$. $B$ induces a holomorphic structure on this bundle and we denote by ${\bL}_\co$ the resulting holomorphic $V$-line bundle.  If we set ${\ve}(\ell)=(1+{\rm sign}\,(\ell))/2$ we then get
\be
SF(\co)=\left\{
\begin{array}{rcl}
 -1-{\ve}(\ell)-2h_0(K_\Sigma -{\bL}_\co)& {\rm if} & n<0 \\
-1-{\ve}(\ell)-2h_0({\bL}_\co) & {\rm if} & n>0 
 \end{array}
 \right.
\label{eq: jet}
\ee

\subsection{Virtual dimensions}
In this final subsection we will  compute virtual dimensions of finite energy moduli spaces.  We will rely heavily on the
techniques of [MMR].

Consider a  $4$ manifold  $\hat{N}$ with a cylindrical end  isometric to
$[0,\infty)\times N$ where $N$ is  disjoint union of  Seifert manifolds
$\{N_j\; ;\; j=1,\ldots, n\}$ of rational degrees $\ell_j$, singularities $(\vec{\alpha}_j, \vec{\beta}_j$ over Riemann surfaces
$\Sigma_j$ of genera $g_j$. To ease the notational burden we will consider only the case $m=1$ (i.e. the boundary is connected). Thus we will drop the index $j$ in the notation of the above objects.  Assume $N$ has $m$ singular fibers with with singularities  ($\vec{\alpha}, \vec{\beta})$.

 Fix a $spin^c$ structure $\hat{\si}$ on $\hat{N}$. This induces a $spin^c$ structure $\sigma$ on $N$ determined by a line bundle $L$. The metric and compatible connections on the end of $\hat{N}$ are prescribed as
indicated in \S 2.2.   This means  that as background connection on $N$ we use
the adiabatic connection $\nabla^\infty$. Consider  a finite energy, irreducible, solution $\hco=(\hat{\phi}, \hat{A})$ of the Seiberg-Witten equations  associated to the structure $\hat{\si}$ and we assume  that along the neck it is in temporal gauge 
\[
\hco=\{t\mapsto\co(t)=(\phi(t),A(t))\}.
\]

The techniques of [MMR] work with no essential changes in the Seiberg-Witten
context and   show that  $[\co(t)]$ converges to $[\co_\infty] \in
\modu_L$, where by $\modu_L$ we denoted the Seiberg-Witten moduli space determined by the $spin^c$ structure $\si$ on $N$. The first conclusion we draw from this fact is that $L$ must be a pulled back  line bundle e since
otherwise $\modu_L =\emptyset$. Suppose this is indeed the case and set $\kappa = c_1(L)$.

The moduli space $\modu_\kappa$ is a disjoint union
\[
\modu_\kappa =\bigcup_n \modu_{\kappa, n}.
\]
Assume first that the configuration $\co_\infty$ is irreducible 
\[
\co_\infty\in \modu_{\kappa,n}(\sigma).
\]
Again, to ease the notation we set $\co =\co_\infty$. For simplicity assume $n<0$ so that is $\co=(\phi_-\oplus, \phi_+,A)$ then $\phi_-=\equiv 0$. The configuration $\co$ is pulled back from a vortex pair $(\tilde{\phi}, B)$ corresponding to a $V$ line bundle $L_\Sigma =L(c,\vec{\gamma}) \in R_\kappa$.

We are interested in  describing a neighborhood of $\hco$ in the moduli space of
finite energy solutions and we will begin as in [MMR] by  studying a simpler
problem. 

Define $\hmodu([\co])$ as the moduli space of  finite energy solutions 
with asymptotic limit $[\co_\infty]$. We want to  understand the structure of a
small neighborhood of $\hco \in \hmodu([\co])$.  More precisely, we would like
to compute the virtual dimension of such a neighborhood. This is achieved  using Kuranishi's deformation picture of the moduli space which  requires  a suitable  functional framework.

Since the convergence to the asymptotic limit  is exponential one can use the
very convenient   weighted Sobolev spaces $L^{k,p}_w$ where $w$ is a very small
positive number.   The resulting  deformation complex can be described as in
Chap.8 of [MMR] and is
\be
0\ra
X_0\stackrel{ \hat{\Lie}_{\hco} }{\longrightarrow} X_1 \stackrel{\underline{sw} }{\longrightarrow}
X_2 \ra 0.
\label{eq: ku1}
\ee
where $X_0$ is the Lie algebra of the group of gauge  transformations on
$\hat{N}$  exponentially  converging to   1 along the neck
\[
X_0=L_w^{3,2}(\ii\Lambda^0T^*\hat{N})
\]
$X_1$ is the tangent space to the space of configurations of the
$4$-dimensional equations
\[
X_1=L_w^{2,2}(\hat{\bS}^+_{\hat{\si}}\oplus\ii\Lambda^1T^*\hat{N})
\]
$X_2$ is the space of ``obstructions''
\[
X_2=L_w^{1,2}(\hat{\bS}^-_{\hat{\si}}\oplus \ii \Lambda^2_+T^*\hat{N})
\]
$\hat{\Lie}=\hat{\Lie}_{\hco}$ is the infinitesimal  gauge group action at
$\hco$ and $\underline{sw}$ is the linearization  at $\hco$ of the $SW$-equations on
$\hat{N}$.

We can now form the  operator
\[
\hat{\cal O}_w:X_1\ra X_2\oplus X_0,\;\;\hat{\cal O}_w =\underline{sw}\oplus
\hat{\Lie}^{\ast_w}
\]
where ${}^{\ast_w}$ denotes the $L^2_w$-adjoint of $\hat{\Lie}$. This is an
elliptic operator and a computation {\em \`{a} la} [MMR] (Chap. 8) shows that along the neck it exponentially approaches 
\[
\hat{\cal O}_w=something \times(\nabla_t-{\cal O}_w) 
\]
where
\[
{\cal O}_w\left[
\begin{array}{c}
{\dps}\oplus{\ida} \\
{\bf i}f 
\end{array}
\right ] = \left[
\begin{array}{clcl}
D_A\dps  &                               &+ &  {\bc}(\ida)\phi-\ii f\phi   \\
                 &-\ii\ast  d\da +\ii df & +&\dta(\phi, \dps) \\
                  &\ii d^*\da -2w\ii f& + & \ii \im \lan \phi, \dps\ran  
\end{array}
\right]
\]
and $[\phi, A]=[\co_\infty]$.  Note that ${\hhes}_1(\co_\infty)={\cal
O}_w\!\mid_{w=0}$. $\hat{\cal O}_w$ is a Fredholm operator and its index (over
${\bR}$) is equal to the virtual dimension  of a small neighborhood of $[\hco] $
in $\hmodu([\co])$. We conclude as in $\S 3.4$ of [N2] that  the index of $\hat{\cal O}_w$ is equal to the Atiyah-Patodi-Singer  index  (${\rm ind}_{APS}$) of $\hat{\cal O}_w$.

Denote by ${\cal A}$ the   anti-selfduality operator on $\hat{N}$
\[
{\cal A}=d_+\oplus d^*:\ii\Omega^1(\hat{N})\ra \ii\Omega^2(\hat{N})\oplus
\ii\Omega^0(\hat{N}).
\]
Using the connection $\hat{A}$  and the $spin^c$ structure $\hsi$ we can form
the Dirac operator 
\[
\hat{D}_{\hat{A}}:\Gamma(\hat{\bS}^+_{\hat{\si}})\ra
\Gamma(\hat{\bS}^-_{\hat{\si}}).
\]
Along the neck  the direct sum ${\cal N}_{\hco}=\hat{D}_{\hat{A}}\oplus {\cal A}$ has the
$APS$ form
\[
{\cal N}=something \times \left(\nabla_t -{\hhes}_0(\co_\infty)\right). 
\]
Using the excision formula in [N2] we deduce
\be
{\rm ind}_{APS}(\hat{\cal O}_w) = {\rm ind}_{APS}({\cal N})-
SF({\hhes}_0\ra {\hhes}_1) - SF({\hhes}_1\ra {\cal O}_w)
\label{eq: vd0}
\ee 
where  $SF(A\ra B)$ denotes $SF(A+t(B-A),\; t\in[0,1])$.  All the indices
and the spectral flows above are   {\em real} quantities.

We now proceed to determine the three terms in the right-hand side of the above
formula.   

Arguing as in the Appendix D of [N2] we conclude that  the third term above vanishes. The second spectral term can be rewritten as
\be
SF({\hhes}_0\ra {\hhes}_1)= SF([\co])
\label{eq: vd1}
\ee
We denote by $\rho_{asd}$ (resp. $\rho_{dir}$) the index densities of ${\cal
A}$ (resp. $\hat{D}_{\hat{A}}$)
\[
\rho_{asd}= -\frac{1}{2}\left({\bf e}(\hat{N})+{\bf L}(\hat{N})\right).
\]
where ${\bf e}(\hat{N})$ and ${\bf L}(\hat{N})$ denote respectively the Euler
and the $L$-genus forms on $\hat{N}$ constructed using the Levi-Civita
connection. Also
\[
\rho_{dir} = 2 \hat{\bf A}(\hat{\nabla}^0) \wedge \exp(\frac{1}{2}c_1(\det \hat{\si}))
\]
where on $\det \hat{\si}$ we used the connection induced by $\hat{A}$.  The
factor 2 appears since we are interested in the {\em real} index of $\hat{D}$. 
The $\hat{\bf A}$-genus form is  computed using the metric compatible  connection $\hat{\nabla}^0$ which along the neck has the product form $dt\otimes \partial_t+\nabla^0$. We denote by $c(\hat{A})$ the  form
\[
\frac{\ii}{2\pi}\hat{F}_{\hat{A}}
\]
on $\hat{N}$ where $\hat{A}$ denotes the induced connection on $\det \hat{\si}$.

On  a $4$-manifold the above equality has a simpler form
\[
\rho_{dir}=\frac{1}{4}(c(\hat{A})^2 -{\bf L}(\hat{\nabla}^0)).
\]
The $\xi$ invariant of ${\hhes}_0$   is  the sum  $\xi({\cal
A}\!\mid_N)+2\xi(D_A)$ (the factor 2 is present  for reality reasons).
\[
\xi({\hhes}_0)=\frac{1}{2}( \dim_{\bf R}\ker {\hhes}_0 -\eta_{sign} +
2\eta(D_A))
\]
where $\eta_{sign}$ denotes the eta invariant of the odd signature operator. We deduce
\[
{\rm ind}_{APS}({\cal N})=\int_N (\rho_{asd}+\rho_{dir})-\xi({\hhes}_0)
\]
\[
=-\frac{1}{2}\int_{\hat{N}} {\bf e} -\frac{1}{2}\left(\int_{\hat{N}} {\bf L}
-\eta_{sign}\right)+\frac{1}{4}\int_{\hat{N}}(c^2(\hat{A})-{\bf L}({\nah}^0))
\]
\[
-\frac{1}{2}\dim_{\bf R}\ker{\hhes}_0 -\eta(D_A).
\]
In [N1] and [MOY] it is shown that the  kernel of ${\hhes}_0$ is isomorphic to
\[
H^0(N, {\bR})\oplus H^1(N, {\bR})\oplus \ker D_A,\;\;\ker D_A\cong H^0({\bL}_\co)+H^0(K_\sigma-{\bL}_\co).
\]
Using the equality (\ref{eq: etadede}) we deduce
\[
{\rm ind}_{APS}(\hat{\cal O}_w)=-(\chi(\hat{N}) +{\rm sign}\,(\hat{N}))/2 +\frac{1}{4}\int_{\hat{N}}(c^2(\hat{A})-{\bf
L}({\nah}^0))
\]
\[
-\dim_{\bf C}\ker D_{A}-\frac{1}{2}(2g+1) -\frac{\ell}{6} +2S (\vec{\beta},\vec{\alpha};\vec{\gamma}) +d(\vec{\beta},\vec{\alpha};\vec{\gamma})
\]
Using (\ref{eq: vd0}) and (\ref{eq: vd1}) we deduce
\[
{\rm ind}\,(\hat{\cal O}_w)= -\frac{1}{2}(\chi(\hat{N}) +{\rm sign}\,(\hat{N}))
+\frac{1}{4}\int_{\hat{N}}(c^2(\hat{A})-{\bf
L}({\nah}^0))
\]
\be
-\left(\dim_{\bf C}\ker D_{A}+SF([\co]) \right)
-\frac{1}{2}(1+ 2g)-\frac{\ell}{6} +S (\vec{\beta},\vec{\alpha};\vec{\gamma}) +d(\vec{\beta},\vec{\alpha};\vec{\gamma}).
\label{eq: vd2}
\ee
This formula can be further simplified using (\ref{eq: jet})  (with $n<0$). We have
\[
\dim_{\bf C}\ker D_A +SF([\co])=h_0({\bL}_\co)-h_0(K_\Sigma-{\bL}_\co)-1-\ve(\ell)
\]
(Riemann-Roch-Kawasaki)
\[
=\deg |{\bL}_\co|  -g -\ve(\ell).
\]
  We can replace  the integral of ${\bf L}(\nah^0)$ with the integral of ${\bf L}(\hat{N})$ plus a correction term given by a  transgression.  This correction term can be computed exactly as in the proof of the second transgression formula of [N2] and we  get (assuming the radius of a regular fiber of $N$ is $1$)
\be
\int_{\hat{N}}{\bf
L}(\nah^0)-\int_{\hat{N}}{\bf
L}(\hat{N})=\frac{2\ell}{3}(\ell^2-\chi)
\label{eq: vd3}
\ee
where $\chi=-\deg K_\Sigma$.    The signature eta invariant of
$N$ was computed in  [O] and is given by
\[
\eta_{sign}=-{\rm sign}\,(\ell) -\frac{2\ell}{3}(\ell^2-\chi)
+\frac{\ell}{3} - 4 S(\vec{\beta},\vec{\alpha})
\]
where, as in Remark \ref{rem: lcdir}
\[
 S(\vec{\beta},\vec{\alpha}) =\sum_{i=1}^m s(\beta_i,\gamma_i;0,0).
\]
We deduce from (\ref{eq: vd3})
\[
\int_{\hat{N}}{\bf L}(\nah^0)+\eta_{sign}-\int_{\hat{N}}{\bf L}(\hat{N})=\ell/3-{\rm sign}\,(\ell) - 4S(\vec{\beta}, \vec{\alpha}).
\]
The term 
\[
\eta_{sign}-\int_{\hat{N}}{\bf L}(\hat{N})
\]
is equal to $-{\rm sign}\,(\hat{N})$ so that we  deduce
\[
\int_{\hat{N}}{\bf L}(\nah^0)={\rm sign}\,(\hat{N}) +\frac{\ell}{3} -{\rm
sign}\,(\ell) - 4S(\vec{\beta}, \vec{\alpha}).
\]
If we use this equality in (\ref{eq: vd2}) we deduce
\[
{\rm ind}_{APS}(\hat{\cal O}_w)= 
\frac{1}{4}\int_{\hat{N}}c^2(\hat{A})-\frac{1}{4}(2\chi(\hat{N})+3 {\rm
sign}\,(\hat{N}) )
\]
\[
+\left(g+{\ve}(\ell)-\deg|{\bL}_\co| \right)
-\frac{1}{2}(2g+1)-\frac{1}{4}(\ell/3 -{\rm sign}\,(\ell)) -\frac{\ell}{6}
\]
\[
+2S(\vec{\beta},\vec{\alpha},\vec{\gamma}) +d(\vec{\beta}, \vec{\alpha},\vec{\gamma})+S(\vec{\beta},\vec{\alpha}).
\]
\[
=d(\hat{\si}) -\frac{1}{2} +{\ve}(\ell)+\frac{ {\rm sign}\,(\ell)}{4}-\frac{\ell}{4}-\deg|{\bL}_\co | 
\]
\be
+2S(\vec{\beta},\vec{\alpha};\vec{\gamma}) +d(\vec{\beta}, \vec{\alpha},\vec{\gamma})+S(\vec{\beta},\vec{\alpha})
\label{eq: vd4}
\ee
where  
\[
d(\hat{\si}):=\frac{1}{4}\int_{\hat{N}}c^2(\hat{A})-\frac{1}{4}(2\chi(\hat{N})+3 {\rm
sign}\,(\hat{N}) ).
\]
Note that $d(\hat{\si})$ is precisely  the expression computing the virtual dimensions of Seiberg-Witten moduli spaces on closed manifolds.

To find the virtual dimension $\dim_v(\hco)$ of a neighborhood of $\hco$ in the {\em entire}
moduli space $\hmodu$  we only have to add the dimensions of the  asymptotic
limit sets $\dim {\modu}_{\kappa, n}=2\deg |{\bL}_\co|$ (recall that we have
assumed $n<0$).
\[
\dim_v(\hco)= d(\hat{\si}) -\frac{1}{2} +{\ve}(\ell)+\frac{ {\rm sign}\,(\ell)}{4}-\frac{\ell}{4}+\deg|{\bL}_\co|
\]
\be
+2S(\vec{\beta},\vec{\alpha},\vec{\gamma}) +d(\vec{\beta}, \vec{\alpha},\vec{\gamma})+S(\vec{\beta},\vec{\alpha})
\label{eq: vd5}
\ee
We can now define the  boundary correction term
\[
\omega([\co]):= -\frac{1}{2} +{\ve}(\ell)-\frac{ \ell -{\rm sign}\,(\ell)}{4}+2S(\vec{\beta},\vec{\alpha},\vec{\gamma}) +d(\vec{\beta}, \vec{\alpha},\vec{\gamma})+S(\vec{\beta},\vec{\alpha})+\deg |{\bL}_\co|.
\]
where $[\co] \in \modu_{\kappa, L}$, where $L\ra \Sigma$ is a $V$-line bundle in $R_\kappa$ with singularity data $\vec{\gamma}$.  The case  when the asymptotic limit $\co$ satisfies $\nu(\co)>0$  (that is $\phi_+\equiv 0)$ can be safely left to the reader.  The only changes in $\omega(\co)$ occur at    the term $\deg |{\bL}_\co|$ above  which should be replaced with $\deg |K-{\bL}_\co|$.

 When the boundary of $\hat{N}$ has several components (all Seifert manifolds of various types)   then the asymptotic limit of a finite energy solution on $\hat{N}$ is  a collection of solutions of the Seiberg-Witten equations on each of the components.  The corresponding correction term is the sum of the correction terms   determined by each component. 

We conclude by considering the case when $\hco$ is an irreducible finite energy solution on $\hat{N}$ with  asymptotic limit $C=C_\infty$ a nondegenerate reducible.  We argue exactly as in [N2].  Set $\rho=\rho(L)$ as in $\S 1.2$ and denote by $L_\Sigma=L(c,\vec{\gamma})$ the canonical representative of $L$.

 The convergence to such a  nondegenerate reducible continues to be exponential and
thus we can use the same   functional framework as above. Assume for simplicity
the boundary has only one component.  We distinguish two cases.

\medskip

\noindent $\bullet$ $\rho(L) \neq 0$. We have to compute
the $APS$ index  of a new operator $\hat{\cal O}_w$ which  along the neck has
the form
\[
\hat{\cal O}_w=something \times (\nabla_t-{\cal O}_w)
\]
where this time
\[
{\cal O}_w\left[
\begin{array}{c}
{\dps}\oplus{\ida} \\
{\bf i}f 
\end{array}
\right ] = \left[
\begin{array}{cl}
D_A\dps  &                                  \\
                 &-\ii\ast  d\da +\ii df  \\
                  &\ii d^*\da -2w\ii f   
\end{array}
\right]
\]
(The spinor part $\phi$ of the asymptotic limit $[\co]=[\phi, A]$ is zero and thus ${\cal P}_\phi\equiv 0$.) Thus  
\[
{\rm ind}_{aps}(\hat{\cal O}_w)= {\rm ind}_{aps}({\cal N})- SF((1-t){\cal O}_0 + t{\cal O}_w).
\]
The spectral flow contribution is easy to determine.     The only eigenvalue of
${\cal O}_{tw}$  contributing to the spectral flow is $-2wt\!\mid_{t=0}$ with
a single eigenfunction $\dps \oplus \ii \da \oplus \ii f$, where $\dps=0$,
$\ii\da=0$ and $f\equiv 1$.  Hence
\[
{\rm ind}_{aps}(\hat{\cal O}_w)= {\rm ind}_{aps}({\cal N}) +1.
\]
The index of ${\cal N}$ can be determined as above using instead the eta invariant of the adiabatic operator coupled with the flat
connection $A \in {\cal R}(L)$ determined in \S 1.2. The nondegeneracy condition also implies $\ker D_A=0$. The eta invariant of $D_A$   is twice its  $\xi$-invariant. Since $\dim_{\bf R} {\hhes}_0= b_0(N)+b_1(N)=1+2g$ we deduce from (\ref{eq: ro1})
\[
{\rm ind}_{APS}(\hat{O}_w)= 1-(\chi(\hat{N}) +{\rm sign}\,(\hat{N}))/2
+\frac{1}{4}\int_{\hat{N}}(c^2(\hat{A})-{\bf
L}({\nah}^0))-\frac{1}{2}(2g+1)
\]
\[
- \frac{1}{2}(\deg K -\deg |K|)(1-2\rho)+\ell\rho(1-\rho) -\frac{\ell}{6} -m\rho
\]
\[
+2S_\rho(\vec{\beta},\vec{\alpha},\vec{\gamma}) +F_\rho(\vec{\alpha},\vec{\beta},\vec{\gamma}).
\]
Again we have
\[
\int_{\hat{N}}{\bf L}({\nah}^0)= {\rm sign}\,(\hat{N}) +\frac{\ell}{3} -{\rm sign}\,(\ell) -4 S(\vec{\beta}, \vec{\alpha}).
\]
We deduce
\[
{\rm ind}_{APS}({\cal O}_w)=d(\hat{\si})+ \frac{1-2g}{2} -\frac{\ell -{\rm sign}\,(\ell)}{4} +\ell\rho(1-\rho)-m\rho
\]
\[
 - \frac{1}{2}(\deg K -\deg |K|)(1-2\rho)+ 2S_\rho(\vec{\beta},\vec{\alpha},\vec{\gamma}) +F_\rho(\vec{\alpha},\vec{\beta},\vec{\gamma}) +S(\vec{\beta},\vec{\alpha}).
\]
 To obtain the virtual dimension of a neighborhood of $\hco$ in the entire moduli space  we need to add the dimension of a neighborhood of ${\co}_\infty$ inside the reducible torus and subtract the dimension of its  isotropy group. This  leads to a correction  by $2g-1$  and   we  deduce
\[
\dim_v(\hco)=d(\hat{\si}) + \omega_{red}(\co)
\]
where
\[
\omega_{red}(\co)=\frac{2g-1}{2} -\frac{\ell -{\rm sign}\,(\ell)}{4}  +\ell\rho(1-\rho)-m\rho- \frac{1}{2}(\deg K -\deg |K|)(1-2\rho) 
\]
\[
 +S(\vec{\beta},\vec{\alpha})+2S_\rho(\vec{\beta},\vec{\alpha},\vec{\gamma}) +F_\rho(\vec{\alpha},\vec{\beta},\vec{\gamma}). 
\]

\noindent $\bullet$ $\rho=0$. In this case we should use the first part of Proposition \ref{prop: eta2}. A computation  as above leads to
\[
\omega_{red}(\co)=\frac{2g-1}{2} -\frac{\ell -{\rm sign}\,(\ell)}{4}  +2S(\vec{\beta},\vec{\alpha},\vec{\gamma})+d(\vec{\beta},\vec{\alpha}; \vec{\gamma}) +S(\vec{\beta},\vec{\alpha}).
\]

\section{Applications}
\setcounter{equation}{0}

\subsection{Tunnelings}
We want to apply the virtual dimension formula formula to the special situation of tunnelings.  Consider a Seifert fibration $N$ with oriented Seifert invariant $(g; \ell \vec{\alpha}, \vec{\beta})$. In this case $\hat{N}={\bR}\times N$ has two boundary components $N_\pm:=\partial_\pm \hat{N}=\{\pm 1\}\times N$ which have opposite orientations. Denote by $(g;\ell_\pm, \vec{\alpha}_\pm,\vec{\beta}_\pm)$ the oriented Seifert invariant of $N_\pm$. We identify $N$ with $N_+$ and thus we have the following identities
\[
\ell_- +\ell_+=0,\;\;\vec{\alpha}_\pm-=\vec{\alpha},\;\; \vec{\beta}_+=\vec{\beta}
\]
and
\[
\vec{\beta}_-=\vec{\alpha}-\vec{\beta}= (\alpha_1-\beta_1, \ldots, \alpha_m-\beta_m).
\]

  Consider a tunneling $\hco$ with irreducible  asymptotic limits $[\co_\pm] \in \modu_{\kappa_\pm, L_\pm}(N)$. Set $n_\pm =\nu(L_\pm)$.  Clearly $\kappa_-=\kappa_+=\kappa$. In [MOY] it is shown that
\[
{\rm sign}\,(n_-)={\rm sign}\,(n_+).
\]
Assume for simplicity that both $n_\pm <0$. Set $k=(\deg(L_+ -\deg L_-)/\ell=(n_+ - n_-)/\ell$. Since $L_\pm \ra \Sigma$ pull back to the same line bundle $\kappa$ we deduce that $k\in {\bZ}$. In fact, $L_+=L_- + k{\bL}_0$ (where $N$ is viewed as the unit sphere bundle of ${\bL}_0$). 

We can now use the  equality (\ref{eq: vd5})   which  gives the virtual dimension  of the space  of tunnelings between $\modu_{\kappa, n_\pm}$
  Denote this dimension by $\tau(\kappa; n_-,n_+)$. We have
\[
\tau(\kappa; n_1,n_2)=\frac{1}{4}\int_{\hat{N}}c(\hat{\si})^2  +\omega([\co_-])+\omega([\co_+]).
\]
The integral term can be computed via transgression exactly as in the third
transgression formula in [N2] and we obtain 
\[
\frac{1}{4}\int_{\hat{N}}c(\hat{\si})^2 =\frac{n_-^2-n_+^2}{\ell}.
\]
Thus we get
\[
\tau(\kappa ;n_-,n_+)=\frac{n_-^2-n_+^2}{\ell}+\omega([\co_-]) +\omega([\co_+]).
\]
The last formula can be made a little bit more transparent. Denote by $\vec{\gamma}_\pm$ the isotropies of $L_\pm$. If
\[
\vec{\gamma}_-=\vec{\gamma}=(\gamma_1, \ldots, \gamma_m)
\]
and $c=\deg L_-$ then $\deg L_+ =c+ k \ell$ and the data $\gamma_+$ are determined by
\[
\vec{\gamma}_+=(\gamma_1^+,\ldots, \gamma_m^+) \;\;{\rm where}\;\; \gamma_i^+=\alpha_i\cdot\{ \frac{\gamma_i+k\beta_i}{\alpha_i}\}.
\]
The correction  term can be further simplified since $\ell_- +\ell_+ =0$, $\ve(\ell_-)+\ve(\ell_+)=1$ and $S(\vec{\beta}_-, \vec{\alpha})+ S(\vec{\beta}_+,\vec{\alpha})=0$.  We conclude
\[
\tau(\kappa ;n_-,n_+)=\frac{n_-^2-n_+^2}{\ell} + \deg |L_-|+\deg|L_+| 
\]
\[
+S(\vec{\beta}_-,\vec{\alpha}, \vec{\gamma}_-) +  S(\vec{\beta}_+,\vec{\alpha}, \vec{\gamma}_+) + d(\vec{\beta}_-,\vec{\alpha}, \vec{\gamma}_-) +  d(\vec{\beta}_+,\vec{\alpha}, \vec{\gamma}_+).
\]
Numerical experiments  show that the above formula is in perfect agreement with the alternative description in [MOY]. 

\medskip

Suppose now that the limit $C_+$ is a nondegenerate reducible. Set as above $n_-=\nu(C_-)$ and denote by $\tau(C_-, C_+)$ the virtual dimension of the space of tunnelings from the component of $C_-$ to the component of $C_+$. We distinguish two cases.

\medskip 

\noindent {\sf A.} $n_- <0$. Arguing as above we deduce 
\[
d(C_-,C_+)=\frac{n_-^2}{\ell} +\ve(-\ell) + g + 2S(\vec{\beta}_-,\vec{\alpha};\vec{\gamma}_-) +d(\vec{\beta}_-,\vec{\alpha};\vec{\gamma}_-) +\deg |{\bL}_{\co_-}| +\delta_+(\rho)
\]
where
\[
\delta_+(\rho)=\ell\rho(1-\rho)- m\rho- \frac{1}{2}(\deg K -\deg |K|)(1-2\rho)
\]
\[
 +2S_\rho(\vec{\beta_+},\vec{\alpha},\vec{\gamma_+}) +F_\rho(\vec{\alpha},\vec{\beta_+},\vec{\gamma}),\;\;{\rm if}\;\rho \neq 0
\]
and 
\[
\delta_+(\rho)= 2S(\vec{\beta_+},\vec{\alpha},\vec{\gamma_+})+d(\vec{\beta_+},\vec{\alpha}; \vec{\gamma_+}), \;\;\;{\rm when}\; \rho=0.
\]

\noindent {\sf B.} $n_->0$.  The  dimension formula is obtained from the above by performing the change
\[
{\bL}_{\co_-}\ra K_\Sigma-{\bL}_{\co_-}.
\]
Of course one should replace $\gamma_+$ with the singularity data of  $K_\Sigma-L$ but as explained in  Remark \ref{rem: duality} this has no effect on the  eta invariants and   thus the overall contributions of the intervening  Dedekind-Rademacher sums is  unchanged by such a substitution.

\subsection{Froyshov invariants of Brieskorn spheres}
 As in [MOY], we want to   consider in some detail the special case of  the Seifert homology spheres $\Sigma(a,b,c)$ where  $a,b,c$ are pairwise coprime  integers $\geq 2$.      The variational theory for the Seiberg-Witten  energy functional  on these manifolds is as simple as one can hope for.  In particular,  using the tunneling informations in the previous section  and our knowledge of eta invariants we will be able to produce  estimates for the Froyshov invariant of such homology spheres.  We begin by  surveying the    geometric facts  which will be relevant in the sequel.

$\Sigma(a,b,c)$ is a Seifert fibration over  an orbi-sphere $S$. We orient it as the link of the complex singularity
\[
\{x^a+y^b+z^c=0\; ;\; |x|^2+|y|^2+|z|^2 \leq {\ve}\}\subset {\bC}^3.
\]
As such, its rational degree is negative $\ell =-\frac{1}{abc}$ and its singularity data are
\[
\vec{\alpha}=(a,b,c),\;\;\vec{\beta}=(\beta_1,\beta_2,\beta_3)
\]
with $\beta_i$ determined  by the congruence
\[
\frac{\beta_i}{\ell\alpha_i} \equiv 1 \;({\rm mod}\; \alpha_i).
\]
The canonical line bundle $K_S$ has rational degree
\[
\kappa = 1-(\frac{1}{a}+\frac{1}{b}+\frac{1}{c}).
\]
Consider the ``simplex''
\[
\Delta (a,b,c) =\{(x,y,z) \in {\bZ}_+^3\; ;\;  \frac{x}{a}+\frac{y}{b}+\frac{z}{c}< \frac{\kappa}{2},\;x<a,\;y<b,\;z<c \}.
\]
To a point ${\bf p}=(x,y,z)\in \Delta(a,b,c)$ we associate a line $V$-bundle $L_{\bf p}\ra S$  with  $\deg |L_{\bf p}| =0$ and  singularity data $(x,y,z)$. Define the ``energy'' of a point ${\bf p}$ by
\[
E({\bf p}) := (\deg L_{\bf p} -\frac{1}{2}\deg K_S)^2/\ell =\nu(L_{\bf p})^2/\ell.
\]

Using the description  of the critical set  of the Seiberg-Witten functional  we deduce that to any ${\bf p}\in \Delta(a,b,c)$  we can associate  a pair  of irreducible solutions $(C_+({\bf p}), C_-({\bf p}))$ of identical energy $E({\bf p})$ (we deliberately omitted the normalization constant in ${\bf Fact 4}$, \S 2.2).  $C_+({\bf p})$ corresponds  to a holomorphic vortex  on $L_{\bf p}$ while $C_-({\bf p})$ corresponds to an antiholomorphic  vortex on $K_S-L_{\bf p}$.    Any irreducible solution belongs to  such a pair. We will denote the unique reducible solution by $C_0$. The discussion at the end of \S 1.2 shows this is nondegenerate.

For each ${\bf p}\in \Delta(a,b,c)$   set
\[
n_\pm({\bf p}) = -\tau(C_\pm({\bf p}), C_0)
\]
where $\tau$ is described in the previous subsection. The dimension formul{\ae} coupled with the Serre duality  for the Dedekind-Rademacher sums (Remark \ref{rem: duality}) show that
\be
n_-({\bf p})=n_+({\bf p})+1.
\label{eq: serre}
\ee
Now  we can form  two Laurent polynomials
\[
P^\pm_{a,b,c}(T) =\sum_{{\bf p}\in \Delta(a,b,c)}T^{n_\pm({\bf p})}.
\]
 Note that according to (\ref{eq: serre})  $P^+_{a,b,c}=P^-_{a,b,c}$. Set $P_{a,b,c}=P^+_{a,b,c}$.  The polynomial$2P_{a,b,c}$ can be  regarded as the Poincar\'{e} polynomial  of the (irreducible) Seiberg-Witten-Floer homology  of $\Sigma(a,b,c)$ (negative gradings are allowed).   Theorem 10.1.1 of [MOY] shows it is an odd polynomial.

We include below  a few examples \[
\begin{array}{ll}
\Sigma(2,3,5)        & P =0        \\
\Sigma(2,3,7)        & P = T^{-1}  \\
\Sigma(2,3,11)       & P =T^{-1}  \\
\Sigma(2,3,13)       & P =T          \\
\Sigma(3,3,17)       & P =T  \\
\Sigma(3,5,7)        & P =T+T^{-1}   \\
\Sigma(3,5,11)       & P =T^5 +T+T^{-1} \\
\Sigma(3,5,13)       & P =T^3+T^5+T^9 \\
\Sigma(5,7,9)        & P =2T+T^3+T^7+T^9+T^{25}
\end{array}
\]
Given a Laurent polynomial 
\[
P(T)=\sum_{n\in {\bf Z}} a_nT^n
\]
 we define its set of {\em gaps} $\Gamma(P)$  by
\[
\Gamma(P) = \{ m \in {\bZ}_+\; ;\; a_{-(2m+1)}=0\}
\]
and set
\[
m(P) =\min \Gamma(P).
\]
For $P=P_{a,b,c}$ the  invariant $m(a,b,c):=m(P_{a,b,c})$ coincides with the integer $m$ defined at the  beginning of  Section 3 in [Fr].  Denote by $Y_{a,b,c}$  the Froyshov invariant of $\Sigma(a,b,c)$ defined in [Fr]. 
 
\begin{theorem}{\rm For any Seifert homology sphere $\Sigma(a,b,c)$ we have the following inequality
\be
Y_{a,b,c} \leq  Z_{a,b,c}:= 8m(a,b,c) + {\bf F}_{a,b,c}
\label{eq: fr3}
\ee
where ${\bf F}_{a,b,c}$ is the invariant ${\bf F}(\Sigma(a,b,c))$ defined in (\ref{eq: fr1}) and (\ref{eq: fr2}) at the end of \S 1.2.}
\label{th: froy}
\end{theorem} 

\noindent{\bf Proof}\hspace{.3cm} We begin by briefly recalling the definition of the invariant $Y$.  We describe only the case of ${\bZ}$-homology spheres.

Suppose $N$ is a ${\bZ}$-homology sphere. Fix an arbitrary metric  $g$ on  $N$.  By suitably perturbing the  Seiberg-Witten equations on $N$ using as usual  a $1$-form $\nu$ on $N$ as perturbation parameter  we can  arrange that the energy functional has only one  nondegenerate reducible solution and the irreducible are isolated nondegenerate.The resulting gradient flow on ${\bR}\times N$ can also be perturbed  to be in a Morse generic situation.

Denote by $S$ the set of irreducible solutions and by $\theta$ the reducible one. For every $\alpha \in S$  denote by $i(\alpha, \theta)$ the dimension of the space of  tunnelings from $\alpha$ to theta (if this space is nonempty). Define  $m$ as the  smallest nonnegative integer such that there are no tunnelings $\alpha \ra \theta$ with $i(\alpha, \theta)=2m+1$.  Denote by $\eta_{dir}$ the eta invariant of the Dirac operator on $N$ (obtained canonically from the Levi-Civita connection) and by $\eta_{sign}$ the eta invariant of the odd signature operator.  Note that these eta invariants  depend on $g$. Set
\[
Y(N,g):= 8m+ 4\eta_{dir}+\eta_{sign}
\]
and  define the Froyshov invariant $Y(N)$ by
\[
Y(N)=\inf \{Y(N,g,\nu)\; ; \; g,\nu\}.
\]
From the definition it is clear that one can obtain upper estimates on $Y(N)$ as soon as one can produce  a concrete example of a metric on $N$ and perturbation $\nu$  satisfying the required nondegeneracy conditions  and moreover  the quantities $m$, $\eta_{sign}$ and $\eta_{dir}$ are explicitly computable.

We present below one such instance when $N=\Sigma(a,b,c)$. As  metric on $N$ we use the deformation of the Thurston metric discussed in Remark \ref{eq: lcdir}. This is obtained by a rescaling along the fibers of the Thurston metric so that the fibers  become very short, of radius $\approx r \ll1$. Denote this metric by $g_r$ and let ${\bS}$ denote the spinor bundle associated to the unique $spin$ structure on $N$.  For any connection $A$ on $\det {\bS}$  we can now construct two Dirac operators: the Levi-Civita induced ${\dir}_{A,r}$ and $D_{A,r}$ induced by the adiabatic connection on $TN$, as in Remark \ref{rem: lcdir}. Correspondingly we get two   Seiberg-Witten equations: $SW^r_{reg}$ formulated in terms of ${\dir}_{A,r}$ and the adiabatic equations $SW^r_{ad}$.  Denote by $S_{reg}(r)$ (resp. $S_{ad}(r)$ the  set of {\em irreducible} solutions of $SW^r_{reg}$ (resp. $SW^r_{ad}$.). $SW_{reg}$ has a unique reducible the obvious one. In [N1] we established the following facts.

\noindent {\bf Fact 1.} $S^r_{ad}$ is independent of $r$. Denote this set by $S_{ad}$. This set is  described at the beginning of this subsection.

\medskip

\noindent {\bf Fact 2.} $S^r_{reg}\ra S_{ad}$ as $r\ra 0$ in any Sobolev norm defined in terms of the metric $g_1$.  Moreover, the rates of convergence can be estimated effectively in terms of  $r$.

Note that the hessians of $SW^r_{ad}$  at solutions $\alpha\in S_{ad}$ are  invertible  but the norms of their inverses depend upon the metric $g_r$.  In [N1] we also established 

\medskip

\noindent {\bf Fact 3.} The eigenvalues of these hessians are bounded away from zero by a constant {\em independent of} $r$.

\medskip

 The last fact implies that  any $\alpha \in S_{ad}$ admits  a neighborhood containing for  each $r\ll 1$  an unique solution of $SW_{r}$ which must be nondegenerate. Thus, for $r\ll 1$ we have a bijection
\be
\phi_r: S^r_{reg}\ra S_{ad}.
\label{eq: map}
\ee
The map $\phi_r$ is also compatible with the gradings.   To show this we need a bit more notations.

For  any $\alpha \in S^r_{reg}$ denote by $H_{reg}(\alpha,r)$ the hessian of $SW^r_{reg}$ at $\alpha$.  Define  $H_{ad}(\beta, r)$ similarly.  Additionally, we have two hessians $H_{reg}(\theta,r)$ and $H_{ad}(\theta,r)$ at the trivial solution $\theta$.

 $\phi_r$ preserves the gradings if
\be
SF(H_{reg}(\alpha,r) \ra H_{reg}(\beta, r))= SF(H_{ad}(\phi_r(\alpha),r)\ra H_{ad}(\phi_r(\beta),r)\,)
\label{eq: grad1}
\ee
and
\be
SF(H_{reg}(\theta,r) \ra H_{reg}(\alpha,r)) =SF(H_{ad}(\theta,r)\ra H_{ad}(\phi_r(\alpha),r) \,).
\label{eq: grad2}
\ee
The equality (\ref{eq: grad1}) follows from the fact that $H_{red}^r(\alpha)$ is very close to $H_{ad}^r(\phi_r(\alpha))$ and both operators are invertible. Hence
\[
SF(H_{red}(\alpha)\ra H_{ad}(\phi_r(\alpha))\,)=0.
\]
To prove (\ref{eq: grad2}) it suffices to show that
\be
SF(H_{red}(\theta)\ra H_{ad}(\theta))=0.
\label{eq: grad3}).
\ee
This spectral flow was analyzed in \S 2.3 of [N2] and it was shown to be zero.

The upshot of the previous discussion is that  in order to compute the ingredient $m$ in the definition of $Y(\Sigma(a,b,c),g_r)$ we may as well use the adiabatic  Seiberg-Witten equations.   We see that it coincides with $m(a,b,c)$. The inequality (\ref{eq: fr3}) is now obvious.  $\Box$

\bigskip

We now associate to any {\em negative definite} ${\bZ}$-quadratic form $q$ an  integer $\Theta(q)$   defined by
\[
\Theta (q) = {\rm rk}\,(q) + \max \{ q(\xi,\xi)\; ;\; \xi \; {\rm characteristic\; vector} \}.
\]
We list below  a few  properties of this invariant.

\medskip

\noindent {\bf P1}  $\Theta (q)$ is divisible by $8$.

This follows from the arithmetic properties of the characteristic vectors described e.g. in [Se].

\medskip

\noindent {\bf P2}  $\Theta(q) \leq {\rm rk}\,(q)$ with equality if and only if $q$ is  even.

Indeed, one has inequality above iff  $0$ is a characteristic vector i.e. $q$ is even.

\medskip

\noindent {\bf P3}    $\Theta (q) \geq 0$ with equality  iff $q$ is diagonal.

This   highly nontrivial result is proven in [E].

\medskip

\noindent {\bf P4} If $q=q_1\oplus q_2$  and $q_2$ is even 
\[
\Theta(q)= \Theta(q_1) +\Theta(q_2)= \Theta(q_1) +{\rm rk}\,(q_2).
\]
This is quite elementary and can be left to the reader.

\medskip

\noindent {\bf P5}  If $M$ is a $4$-manifold with  boundary a  ${\bZ}$-homology sphere $N$ and if the intersection form $q$ of $M$ is negative definite then
\[
\Theta(q) \leq Y(N)
\]
where $Y(N)$ denotes the Froyshov invariant of $N$.

This is  Theorem  1 in  [Fr].

For every  triple $(a,b,c)$ of pairwise coprime  integers $\geq 2$ denote by $\Gamma_{a,b,c}$ the intersection form of the Hirzebr\"{u}ch-Jung plumbing (see [HNK] or [NR] for a definition) associated to the Brieskorn homology sphere $\Sigma(a,b,c)$. $\Gamma_{a,b,c}$ is negative definite and set 
\[
\Theta_{a,b,c}:= \Theta(\Gamma_{a,b,c}).
\]
Note that Theorem \ref{th: froy} and {\bf P5} imply the  following inequality
\be
0\leq \Theta_{a,b,c} \leq Y_{a,b,c} \leq Z_{a,b,c}
\label{eq: bound}
\ee
The following topological result  is now obvious.

\begin{proposition}{\rm   Suppose $\Sigma(a,b,c)$ is such that $Z_{a,b,c}\leq 0$. Then
\[
Y_{a,b,c}=Z_{a,b,c}=0.
\]
Moreover,  if $X$ is a  $4$-manifold with negative  definite intersection form $q$  and $\partial X= \Sigma(a,b,c)$ then  $q$ is diagonalizable.}
\label{prop: diagonal}
\end{proposition}

Naturally, one may ask whether there exist $\Sigma(a,b,c)'s$ with $Z=0$. We present below a few values of $Z$. 

\medskip
\begin{center}
\begin{tabular}{||c|c|c|c||} \hline
$(a,b,c)$    & {\bf F} &  $8m$   & $Z$   \\ \hline\hline

$(2,3,5)$    &  $8$    &  $0$    & $8$     \\ \hline
$(2,3,7)$    & $-8$    &  $8$    & $0$      \\ \hline
$(2,3,11)$   & $0$     &  $8$    & $8$      \\ \hline
$(2,3,13)$   & $0$     &  $0$    & $0$      \\ \hline
$(2,3,17)$   & $8$     &  $0$    & $8$  \\  \hline\hline 
$(3,5,7)$    & $0$     &  $8$    & $8$  \\ \hline   
$(3,5,11)$   & $0$     &  $8$    & $8$   \\ \hline
$(3,5,13)$   & $8$     &  $0$    & $8$   \\ \hline\hline
$(5,7,9)$    & $0$     &  $0$    & $0$  \\ \hline
\end{tabular}
\end{center}

The periodicity displayed in the first part of the table  is a manifestation  of a more general result.

\begin{proposition}{\rm Denote by $Q_k$, ${\bf F}_k$ and respectively $Z_k$  the quantities $P_{2,3,6k+1}$, ${\bf F}_{2,3,6k+1}$ and respectively $Z_{2,3,6k+1}$. Then
\be
Q_k=\left\{
\begin{array}{rll}
jT^{-1} & ,& k= 2j-1 \\
jT & ,& k=2j
\end{array}
\right.
\label{eq: map1}
\ee
\be 
{\bf F}_k=\left\{
\begin{array}{rll}
-8 &,& k=2j-1 \\
0 &,& k=2j
\end{array}
\right.
\label{eq: map2}
\ee
In particular, $Z_k=0$ for all $k\geq 1$.}
\label{prop: map}
\end{proposition}

\noindent {\bf Sketch of proof} \hspace{.3cm}  Let us first point out that a weaker version of (\ref{eq: map1}) is established in [MOY].   Computing  the Dedekind sums may  in general be a very daunting task.   The    singularity data of $\Sigma(2,3, 6k+1)$ are
\[
\vec{\alpha}=(2,3,6k+1)\;\;{\rm and}\;\;\vec{\beta}=(1,1, k).
\]
The ``simplex'' $\Delta(2,3,6k+1)$ is easily described since
\[
\deg K_S = -\frac{5}{6} +\frac{6k}{6k+1}
\]
is extremely small.   The points in $\Delta(2,3, 6k+1)$ have the form $(0,0,z)$ where $z$ is a nonnegative integer such that $z/(6k+1) < \deg K_S/2$. Moreover, the invariant $\rho$ of the unique line bundle over $\Sigma(2,3,6k+1)$ is $1/2$.

\begin{figure}
\centerline{\psfig{figure=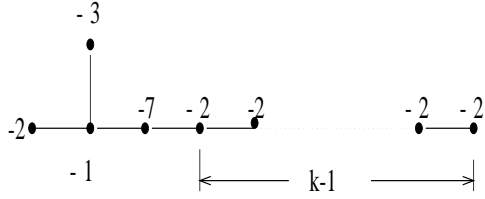,height= 1in,width=2.5in}}
\caption{\sl {The plumbing diagram of $\Sigma(2,3,6k+1)$}}
\label{fig: seif1}
\end{figure}

The computation of the quantities $n_\pm({\bf p})$ and  ${\bf F}$ reduces essentially to  computing  Dedekind-Rademacher sums of the form
\[
s(k, 6k+1, r, 0)\;\;{\rm and}\;\; s(k, 6k+1, r, -1/2)
\]
where $r$ is a rational number.  At this point the reciprocity law comes in very handy. Denote by $R$ the reciprocity function defined by (\ref{eq: rec1}) and (\ref{eq: rec2}). Then
\[
s(k,6k+1;x,y)=-s(6k+1,k;y,x) +R(k,6k+1;x,y)
\]
\[
=-s(1,k;y+6x,x) +R(k,6k+1; x,y)
\]
\[
 =s(k,1;x,y+6x) -R(1,k;,y+6x,x)+R(k,6k+1;x,y)
\]
\[
= ((ky+6kx))\cdot ((y+6x)) -R(1,k;,y+6x,x)+R(k,6k+1;x,y).
\]
Above we used the usual continuous fraction decomposition of $k/(6k+1)$ which, from a computational point of view, is more efficient  than the Hirzebr\"{u}ch-Jung   continuous fraction decomposition. In this case the latter requires $k$ steps (see Figure \ref{fig: seif1}). Obtaining the  identities (\ref{eq: map1}) and (\ref{eq: map2})  is now an elementary, albeit tedious,  accounting job.  It can be safely left to the reader.  $\Box$

\medskip

From  Proposition \ref{prop: map}  and  \ref{prop: diagonal} we deduce immediately the following  consequence.

\begin{corollary}{\rm If $X$ is a smooth $4$-manifold with negative intersection form bounding  $\Sigma(2,3,6k+1)$ then its intersection form must be diagonalizable.}
\label{cor: liviu}
\end{corollary}

\begin{remark}{\rm The above  corollary implies  that $\Gamma_{2,3, 6k+1}$ is diagonalizable.  This consequence is not entirely obvious but  it can be proved  directly.  We present below a simple  argument discovered in conversations with Ian Hambleton.  

Set for simplicity ${\cal I}_k=\Gamma_{2,3,6k+1}$.  More explicitly, ${\cal I}_k$ is given by the plumbing  diagram in Figure \ref{fig: seif1}.   Note that ${\rm rk}\,({\cal I}_k)=4+k-1$.  ${\cal I}_1\cong 4\lan- {\bf 1}\ran$  since all odd, unimodular, negative definite forms of rank $4$ are so. Its matrix is
\[
A=\left[
\begin{array}{cccc}
-1 & 1 & 1 & 1 \\
1 & -2 & 0 & 0 \\
1 & 0 & -3 & 0 \\
1 & 0 & 0 & -7 
\end{array}
\right]
\]
with inverse
\[
B = - \left[
\begin{array}{cccc}
42 & 21 & 14 & 6 \\
21 & 11&  7 & 3 \\
14 & 7 & 5 & 2 \\
6 & 3 & 2 & 1 
\end{array}
\right].
\]
The matrix of  ${\cal I}_k$ can be written in block form
\[
A_k=\left[
\begin{array}{cc}
A & L^t \\
L & N
\end{array}
\right]
\]
where
\[
L=\left[
\begin{array}{cccc}
0 & 0 & 0 & 1\\
0 & 0 & 0 & 0 \\
0 & 0 & 0 & 0 \\
\vdots & \vdots & \vdots & \vdots 
\end{array}
\right]
\]
and
\[
N=\left[\begin{array}{cccccc}
-2 & 1 & 0 & 0  & 0 &\cdots \\
1 & -2 & 1 & 0 & 0 &\cdots \\
0 & 1 & -2 & 1 & 0 & \cdots \\
 &   &  &  &  &  \\

\vdots & \vdots & \vdots & \vdots & \vdots & \vdots
\end{array}
\right]
\]
According to [HNK], Lemma 1.4,   the matrix $A_k$ is equivalent (as  ${\bZ}$-quadratic forms) to
\[
C=    4\lan -{\bf 1} \ran \oplus   (-LBL^t+N).
\]
A simple computation shows that
\[
LBL^t= \left[
\begin{array}{cccc}
-1 & 0  & 0 & \cdots \\
0 & 0 &  0 & \cdots \\
 & & & \\
\vdots & \vdots & \vdots & \vdots 
\end{array}
\right]
\]
Thus $-LBL^t +N$ is the   matrix of the intersection form $\lan -1|k-2\ran$ given by the diagram
\[
\stackrel{-1}{\bullet}-\underbrace{\stackrel{-2}{\bullet}- \; \cdots \; -\stackrel{-2}{\bullet}-\stackrel{-2}{\bullet}}_{k-2}.
\]
Using the same trick again, we can split off  the isolated $-1$ in the above diagram  and one can  verify immediately the  identity
\[
\lan -1 | m\ran \cong  \lan -{\bf 1} \ran \oplus \lan -1 | m-1 \ran.
\]
This shows inductively that the  form ${\cal I}_k$ is diagonal.}
\label{rem: funny}
\end{remark}

\begin{remark}{\rm As pointed out to the author  by S. Akbulut and R. Fintushel, the Brieskorn spheres $\Sigma(2,3,13)$ and $\Sigma(2,3,25)$ bound  contractible smooth manifolds and in these cases Corollary \ref{cor: liviu} follows from  Donaldson's first theorem ([DK], Thm. 1.3.1). Also, as indicated by Nikolai Saveliev, $\Sigma(2,3,7)$ bounds a rational ball (with $2$-torsion in its first homology) and in this case the Corollary \ref{cor: liviu} also follows from Donaldson's first theorem.}
\label{rem: akfin}
\end{remark}

The  properties {\bf P1}-{\bf P5}  and Theorem \ref{th: froy} have another interesting topological consequence.

\begin{proposition}{\rm  Suppose $Z_{a,b,c}\leq 8$ and $X$ is a $4$-manifold  with the following properties.

\noindent (i) $\partial X =\Sigma(a,b,c)$

\noindent (ii) The intersection form $q$ of $X$ is negative definite and splits as $q=q_1\oplus q_2$ where $q_2\neq 0$ is even.

 Then $q_2=-E_8$, $q_1$ is diagonal and $Y_{a,b,c}=8$.}
\label{prop: 8}
\end{proposition}

\noindent {\bf Warning}  The above proposition  does not prove $Z_{a,b,c}=8\Rightarrow Y_{a,b,c}=8$. It only shows this is the cases provided $\Sigma(a,b,c)$ bounds a manifold satisfying (ii) above.

\medskip

\noindent{\bf Proof} \hspace{.3cm} Indeed, from {\bf P4} we deduce
\[
\Theta(q_1)+{\rm rk}\,(q_2) \leq 8.
\]
By {\bf P3}  $\Theta(q_1) \geq 0$  so that ${\rm rk}\,(q_2) \leq 8$ which forces $q_2=-E_8$.  In particular $\Theta(q_1) =0$ so that $q_1$ must be diagonal according to {\bf P3}. Hence $\Theta(q)=8$ Using {\bf P5}  and (\ref{eq: bound}) we deduce
\[
8=\Theta (q) \leq Y_{a,b,c} \leq Z_{a,b,c} \leq 8. \;\;\Box
\]

\medskip

We can  also prove a counterpart of   Proposition \ref{prop: map}.

\begin{proposition}{\rm Denote by $Q_k$, ${\bf F}_k$ and respectively $Z_k$ the quantities $P_{2,3,6k-1}$, ${\bf F}_{2,3,6k-1}$ and respectively $Z_{2,3,6k-1}$. Then
\be
Q_k=\left\{
\begin{array}{rll}
jT^{-1} & ,& k= 2j \\
jT & ,& k=2j+1
\end{array}
\right.
\label{eq: map3}
\ee
\be 
{\bf F}_k=\left\{
\begin{array}{rll}
8 &,& k=2j+1 \\
0 &,& k=2j
\end{array}
\right.
\label{eq: map4}
\ee
In particular, $Z_k=8$ for all $k\geq 1$.}
\label{prop: map2}
\end{proposition}

\noindent{\bf Proof} \hspace{.3cm}   Follows the same lines  as the proof of Proposition \ref{prop: map}. The only difference now   is in the  Seifert invariants of $\Sigma(2,3,6k-1)$. They are
\[
\vec{\alpha}=(2,3,6k-1),\;\;\vec{\beta}=(1,2,5k-1).
\]
The Dedekind-Rademacher reciprocity  law drastically reduces the computations    using the following (unorthodox) continuous fraction decomposition of $(5k-1)/(6k-1)$
\[
\frac{5k-1}{6k-1}=\frac{1}{ 1 +\frac{1}{5 -\frac{1}{k}}}. \;\;\Box
\]

\medskip

\begin{corollary}{\rm  (i) $Y_{2,3,6k-1}=8$ for all $k\geq 1$.

\noindent (ii) If  $\partial X =\Sigma(2,3,6k-1)$ and the intersection form $q$ of $X$ is negative definite and splits as $q=q_1\oplus q_2$, $q_2\neq 0$ even, then $q_1$ is diagonalizable and $q_2\cong -E_8$.}
\label{cor: liviu1}
\end{corollary}

\noindent{\bf Proof}\hspace{.3cm} To prove part (i) we will show that the Hirzebr\"{u}ch-Jung plumbing $X_{2,3,6k-1}$ satisfies all the conditions in  Proposition \ref{prop: 8}.   Since $Z_{2,3,6k-1}=8$ this will prove the equality $Y_{2,3,6k-1}=8$.

\begin{figure}
\centerline{\psfig{figure=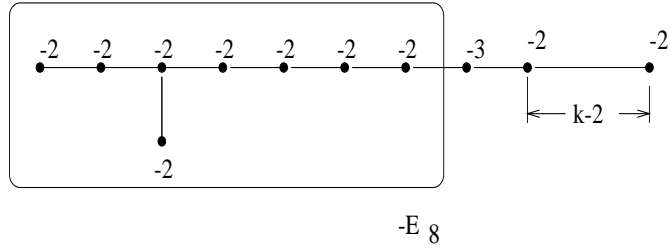,height= 1.3in,width=3.5in}}
\caption{\sl {The plumbing diagram of $\Sigma(2,3,6k-1)$}}
\label{fig: seif2}
\end{figure}

The intersection form $\Gamma_{2,3,6k-1}$  can be read from the plumbing diagram in Figure  \ref{fig: seif2}.   Arguing  as in Remark \ref{rem: funny} we deduce that $\Gamma_{2,3,6k-1} ={\rm diagonal} \oplus (-E_8)$.   Thus $Y_{2,3,6k-1}=Z_{2,3,6k-1}=8$.  The second part of  Corollary \ref{cor: liviu1}  is a special case of Proposition \ref{prop: 8}.  $\Box$

\begin{remark}{\rm The case $k=1$ in the above corollary was  first proved in [Fr].}
\end{remark}

\begin{remark}{\rm Set
\[
d_n= 0+4+8+\cdots 4(n-1)=4n(n-1)
\]
and define 
\[ 
D_n(T)=\sum_{k=1}^n T^{d_k}.
\]
Numerical experimentations reveal  a  beautiful structure of the  polynomials $P_k:=P^+_{2,4k+1,4k+3}$. More precisely they show
\[
P_1(T):=D_1(T)=T, \;\;P_{k}(T)=P_{k-1}+TD_k(T).
\]
Here are the first few  of the polynomials $T^{-1}P_k$
\[
\begin{array}{c}
1\\
T^4+2 \\
T^{12}+2T^4+3 \\
T^{24}+2T^{12}+3T^4+4\\
T^{40}+2T^{24}+3T^{12}+ 4T^4+5 \\
\vdots
\end{array}
\]
Numerical experimentations also show that ${\bf F}_{2,4k+1,4k+3}=0$ so that
$Z_{2, 4k+1,4k+3}=Y_{2,4k+1,4k+3}=0$.   Thus, the only negative definite $4$-manifolds which can bound $\Sigma(2,4k+1,4k+3)$ must have diagonalizable  intersection form.  One can give an alternative proof of this fact. More precisely, according to [CH] the Brieskorn manifold $\Sigma(2,4k+1,4k+3)$ bound s a contractible manifold $B(2,4k+1,4k+3)$.   If $X$ is negative definite  and bounds $\Sigma(2,4k+1,4k+3)$ then
\[
M= X \# B(2,4k+1,4k+3)
\]
is closed and negative definite and by Donaldson's theorem.  A strikingly similar pattern is observed for the polynomials $P^+_{3,3s+1,3s+2}$.  We list the first few of them ($s=1,...,4$) and let the reader deduce the recurrence
\[
\begin{array}{c}
T\\
T^7+2T \\
T^{19}+2T^7+3T \\
T^{37}+2T^{19}+3T^7+4T\\
\vdots 
\end{array}
\]
Again ${\bf F}_{3,3k+1,3k+2}=0$ so that $Y_{3,3k+1,3k+2}=0$. This agrees again with  the conclusions of [CH] where it is established the spheres $\Sigma(3,3k+1,3k+2)$ bound contractible  manifolds.}
\label{rem: 0}
\end{remark}

\subsection{Concluding remarks}

The results we have established so far raise the following natural questions.

\medskip

\noindent {\em   Is it true that $Y_{a,b,c}=Z_{a,b,c}=\Theta_{a,b,c}$ for all $(a,b,c)$? }

The inequality (\ref{eq: bound})   shows the answer is positive provided one can establish only the purely number theoretic identity
\[
Z_{a,b,c}=\Theta_{a,b,c}.
\]
To answer this question  a more manageable  description of $Z_{a,b,c}$ and $\Theta_{a,b,c}$ is highly desirable.

\medskip

Remark \ref{rem: 0} suggests the following question

\noindent{\em Is it true that any $\Sigma(a,b,c)$ which bounds a contractible manifold must have $Y_{a,b,c}=0$?}

Another number theoretic  problem  with possible interesting topological  implications is the following.

{\em Is it true  that for any $\kappa \geq 0$ there exist only finitely many nonisomorphic, irreducible, negative definite, ${\bZ}$-quadratic forms $q$ satisfying}
\[
\Theta(q) \leq \kappa?
\]
In [E] it is proven the answer is positive if $\kappa <16$.

The Euler characteristic of the Seiberg-Witten-Floer homology  of $\Sigma(a,b,c)$ (with the Thurston metric) is $2P_{a,b,c}(-1)$. Following [C1-2] and [KM] define 
\[
\alpha(a,b,c)=2P_{a,b,c}(-1)-\frac{1}{8}{\bf F}_{a,b,c}.
\]
(The  negative sign in front of ${\bf F}_{a,b,c}$  follows from   our ``orientation'' of the Dirac operator  which, although is the  canonical one in the sense of [BC],  it is not the one traditionally used in Seiberg-Witten-Floer theory.) Note that
\[
\alpha(2,3,6k+1) = -k =\alpha(2,3,6k-1)=-k.
\]
The calculations in [FSt]  imply that
\[
\alpha(2,3,6k\pm 1) =\frac{1}{8}\sigma(2,3,6k\pm 1)=
\lambda(2,3,6k\pm 1)
\]
where  $\sigma(a,b,c)$ denotes the signature of the Milnor fiber of $\Sigma(a,b,c)$ and  $\lambda(a,b,c)$ denotes the Casson invariant of $\Sigma(a,b,c)$.

We now define for any oriented homology $3$-sphere $N$ and any  good metric $g$ on $N$ (in the sense of [Mar]) the  quantity
\[
\alpha (N,g)= \chi(N,g) -\frac{1}{8}{\bf F}(N,g)
\]
where $\chi(N,g)$ denotes the Euler characteristic of the  Seiberg-Witten Floer homology  of $(N,g)$ (defined in [Mar]) and 
\[
{\bf F}(N,g)=4\eta_{dir}(g) +\eta_{sign}(g).
\]
The  wall crossing formula of [Mar]  and the variational formul{\ae} for the eta invariants  show  that $\alpha(N,g)$ is independent of the metric  and thus it is a topological invariant of  $N$. 

In [KM], P.Kronheimer and T. Mrowka   conjectured  that this invariant coincides with the Casson invariant and  announced a  proof of its validity when $N$ is a Brieskorn 3-sphere.    In this case the conjecture is equivalent to an interesting  number theoretic identity
\[
\frac{1}{8}\sigma(a,b,c)+2P_{a,b,c}(1)=-\frac{1}{8}{\bf F}_{a,b,c}.
\]
 The right hand side can be expressed in terms    numbers of lattice points inside some  convex regions of ${\bZ}^3$ (see [FSt], [HZ])  while the right hand side  is expressed in terms of Dedekind-Rademacher sums.

\end{document}